\documentclass[twoside,draft,12pt]{article}
\usepackage{amssymb}

\usepackage{makeidx}
\usepackage{a4wide}
\usepackage{amsthm}
\usepackage{amsmath}
\usepackage[all]{xy}
\usepackage{enumerate}
\usepackage{fancyhdr}

\newtheorem{theorem}{Theorem}[section]
\newtheorem{proposition}[theorem]{Proposition}
\newtheorem{corollary}[theorem]{Corollary}
\newtheorem{lemma}[theorem]{Lemma}
\theoremstyle{definition}
\newtheorem{example}[theorem]{Example}
\newtheorem{c-example}[theorem]{Counter Example}
\newtheorem*{Notation}{Notation}
\newtheorem*{Beweis}{Proof}
\newtheorem{definition}[theorem]{Definition}
\newtheorem{punto}[theorem]{}
\theoremstyle{remark}
\newtheorem{remark}[theorem]{Remark}
\newtheorem{remarks}[theorem]{Remarks}
\CompileMatrices
\input{tcilatex}

\begin{document}

\title{Hopf Pairings and (Co)induction Functors over Commutative Rings\thanks{%
MSC (2000): 16W30, 14L17, 16W35, 20G42 \newline
Keywords: Hopf Pairings, Bialgerba Pairings, Dual Bialgebras, Dual Hopf
Algebras, Affine Group Schemes, Quantum Groups, Coinduction Functors,
Cotensor Product.}}
\author{\textbf{Jawad Y. Abuhlail}\thanks{%
Current Address: Box $\#$ 281 King Fahd University of Petroleum $\&\;$%
Minerals, 31261 Dhahran (Saudi Arabia), \textbf{email:}
abuhlail@kfupm.edu.sa, \textbf{Homepage:}
http://faculty.kfupm.edu.sa/math/abuhlail} \\
Mathematics Department\\
Birzeit University\\
P.O.Box 14, Birzeit - Palestine}
\date{}
\maketitle

\begin{abstract}
(\emph{Co})\emph{induction functors} appear in several areas of Algebra in
different forms. Interesting examples are the so called induction functors
in the Theory of Affine Algebraic Groups. In this paper we investigate Hopf
pairings (bialgebra pairings) and use them to study (co)induction functors
for affine group schemes over arbitrary commutative ground rings. We present
also a special type of Hopf pairings (bialgebra pairings) satisfying the so
called $\alpha $-condition. For those pairings the \emph{coinduction functor}
is studied and nice descriptions of it are obtained. Along the paper several
interesting results are generalized from the case of base fields to the case
of arbitrary commutative (Noetherian) ground rings.
\end{abstract}

\section*{Introduction}

\emph{Hopf pairings} (respectively \emph{bialgebra pairings})\emph{\ }were
presented by M. Takeuchi \cite[Page 15]{Tak77} (respectively S. Majid
\cite[1.4]{Maj90}). With the help of these, several authors studied affine
group schemes and quantum groups over arbitrary commutative ground rings
(e.g. \cite{FRT89}, \cite{Tak92}, \cite{Tak95}). In this paper we study the
category of Hopf pairings $\mathcal{P}_{Hopf}$ and the category of bialgebra
pairings{\normalsize \ $\mathcal{P}_{Big}$} over an arbitrary commutative
base ring. In the case of Noetherian base rings we present the \emph{full }%
subcategories $\mathcal{P}_{Hopf}^{\alpha }\subset \mathcal{P}_{Hopf}$
(respectively $\mathcal{P}_{Big}^{\alpha }\subset \mathcal{P}_{Big}$) of
Hopf pairings (respectively bialgebra pairings) satisfying the so called $%
\alpha $\emph{-condition, }see \ref{alpha}. For those a \emph{coinduction
functor} is presented and an interesting description of it is obtained.

The paper is divided into seven sections. The first section includes some
preliminaries about the so called \emph{measuring }$\alpha $\emph{-pairings}%
, \emph{rational modules} and \emph{dual coalgebras}. In the
second section we consider the cotensor functor and prove some
properties of it that will be used in later sections. In the third
section we introduce the coinduction functor in the category of
measuring $\alpha $-pairings and prove mainly that it can be
obtained as a special case of a more general coinduction functor
between categories of type $\sigma \lbrack M].$ Hopf pairings
(bialgebra pairings) are presented in the fourth section, where an
algebraically topological approach is used and several duality
theorems are
proved. In the fifth section we consider the category of Hopf $\alpha $%
-pairings (bialgebra $\alpha $-pairings) and generalize known results on%
\emph{\ admissible }Hopf pairings over Dedekind domains to the case of \emph{%
quasi admissible }bialgebra $\alpha $-pairings and Hopf $\alpha $-pairings
over arbitrary commutative Noetherian ground rings. There the coinduction
functor is also studied and different forms of it that appear in the
literature are shown to be equivalent. The classical duality between groups
and commutative Hopf algebras (e.g. \cite{Mon93}, \cite{Swe69}) is the
subject of the sixth section. In the seventh and last section we apply
results obtained in the previous sections to affine group schemes over
arbitrary commutative rings.

Throughout this paper $R$ denotes a commutative ring with $1_{R}\neq 0_{R}.$
We consider $R$ as a left (and a right) linear topological ring with the
discrete topology. All $R$-modules are assumed to be unital and category of $%
R$-(bi)modules will be denoted by $\mathcal{M}_{R}.$ With unadorned $\mathrm{%
Hom}(-,)$ and $-\otimes -$ we mean $\mathrm{Hom}_{R}(-,)$ and $-\otimes _{R}-
$ respectively. For an $R$-module $M$ we call an $R$-submodule $K\subset M$
\emph{pure} (in the sense of Cohn), if the canonical mapping $\iota
_{K}\otimes _{R}id_{N}:K\otimes _{R}N\rightarrow M\otimes _{R}N$ is
injective for every $R$-module $N.$ For an $R$-module $M$ and subsets $%
X\subset M$ (respectively $Y\subset M^{\ast }$) we set
\begin{equation*}
\mathrm{An}(X):=\{f\in M^{\ast }:f(X)=0\}\text{ (respectively }\mathrm{Ke}%
(Y):=\{m\in M:f(m)=0\text{ }\forall \text{ }f\in Y\}\text{).}
\end{equation*}

Let $A$ be an $R$-algebra (not necessarily with unity). A left $A$-module is
said to be \emph{faithful}, if $\mathrm{An}(_{A}M):=\{a\in A:aM=0\}=(0_{A}).$
We define a left $A$-module $M$ to be $A$\emph{-faithful }(respectively
\emph{unital}), if the canonical map $\rho :M\rightarrow \mathrm{Hom}%
_{R}(A,M)$ is injective (respectively, if $AM=M$). With $\widetilde{_{A}%
\mathcal{M}}$ (respectively $_{A}\mathcal{M}$) we denote the category of $A$%
-faithful (respectively \emph{unital}) left $A$-modules and left $A$-linear
maps. The categories of\emph{\ }$A$-faithful (respectively unital) right $A$%
-modules $\widetilde{\mathcal{M}_{A}}$ (respectively $\mathcal{M}_{A}$) are
analogously defined. If $A$ has unity, then obviously every unital left (or
right) $A$-module is $A$-faithful. For an $R$-algebra $A$ and an $A$-module $%
M,$ an $A$-submodule $N\subset M$ will be called $R$\emph{-cofinite}, if $M/N
$ is finitely generated in $\mathcal{M}_{R}.$ Unless otherwise explicitly
mentioned, we assume that all $R$-algebras have unities respected by $R$%
-algebra morphisms and that all modules of $R$-algebras are unital.

We assume the reader is familiar with the theory and notation of Hopf
Algebras. For any needed definitions the reader may refer to any of the
classical books on the subject (e.g. \cite{Abe80}, \cite{Mon93}, \cite{Swe69}%
) or to the recent monograph \cite{BW03} for the theory of coalgebras over
arbitrary base rings. For an $R$-coalgebra $C,$ we call a right
(respectively a left) $C$-comodule $(M,\varrho _{M})$ \emph{counital }if its
structure map $\varrho _{M}$ is injective, compare \cite[Lemma 1.1.]{CC94}.
For an $R$-coalgebra $C$ we denote with $\mathcal{M}^{C}$ (respectively $^{C}%
\mathcal{M}$) the category of \emph{counital }right (respectively left) $C$%
-comodules. For an $R$-coalgebra $(C,\Delta _{C},\varepsilon _{C})$ and an $R
$-algebra $(A,\mu _{A},\eta _{A})$ we consider $\mathrm{Hom}_{R}(C,A)\;$as
an $R$-algebra with multiplication the so called \emph{convolution product }$%
(f\star g)(c):=\sum f(c_{1})g(c_{2})$ and unity $\eta _{A}\circ \varepsilon
_{C}.$

\section{Preliminaries}

In this section we present some definitions and lemmas to be referred to
later in the paper.

\begin{punto}
\textbf{Subgenerators}.\label{subg} Let $A$ be an $R$-algebra (not
necessarily with unity) and $K$ be a left $A$-module. We say a left $A$%
-module $N$ is $K$\emph{-subgenerated}, if $N$ is isomorphic to a submodule
of a $K$-generated left $A$-module (equivalently, if $N$ is kernel of a
morphism between $K$-generated left $A$-modules). The \emph{full }%
subcategory of $_{A}\mathcal{M},$ whose objects are the $K$-subgenerated
left $A$-modules is denoted by $\sigma \lbrack _{A}K].$ In fact $\sigma
\lbrack _{A}K]\subseteq $ $_{A}\mathcal{M}$ is the \emph{smallest }%
Grothendieck full subcategory that contains $K.$ If $M$ is a left $A$%
-module, then
\begin{equation*}
\mathrm{Sp}(\sigma \lbrack _{A}K],M):=\sum \{f(N)\mid \text{ }f\in \mathrm{%
Hom}_{A-}(N,M),\text{ }N\in \sigma \lbrack _{A}K]\}
\end{equation*}
is the largest $A$-submodule of $M$ that belongs to $\sigma \lbrack _{A}K].$
The subcategory $\sigma \lbrack _{A}K]\subseteq $ $_{A}\mathcal{M}$ can also
be seen as the category of discrete left $A$-modules, where $A$ is
considered as a left linear topological $R$-algebra with the $K$\emph{-adic
topology} (e.g. \cite{Ber94}). The reader is referred to \cite{Wis88} and
\cite{Wis96} for the well developed theory of categories of this type.
\end{punto}

\qquad An important result to which we will often refer is

\begin{lemma}
\label{dense}\emph{(\cite[15.8]{Wis88}, \cite[42.2]{BW03})} Let $A$ be a
ring, $K$ be a faithful left $A$-module and $B\subset A$ be a subring. Then $%
\sigma \lbrack _{B}K]=\sigma \lbrack _{A}K]$ if and only if $B\subset A$ is $%
K$-dense.
\end{lemma}

\begin{remark}
Let $C$ be an $R$-algebra. Then $C^{\ast }$ becomes two (left) linear
topologies, the $C$-adic topology induced by $_{C^{\ast }}C$ and the finite
topology induced by the embedding $C^{\ast }\hookrightarrow R^{C}.$ By
\cite[Lemma 2.2.4]{Abu2001} the two topologies coincide.
\end{remark}

\begin{punto}
\textbf{The }$\alpha $\textbf{-condition}.\label{alpha} With an $R$\emph{%
-pairing} $P=(V,W)$ we mean $R$-modules $V,W$ with an $R$-linear map $\kappa
_{P}:V\rightarrow W^{\ast }$ (equivalently $\chi _{P}:W\rightarrow V^{\ast }$%
). For $R$-pairings $(V,W)$ and $(V^{\prime },W^{\prime })$ a morphism $(\xi
,\theta ):(V^{\prime },W^{\prime })\rightarrow (V,W)$ consists of $R$-linear
mappings $\xi :V\rightarrow V^{\prime }$ and $\theta :W^{\prime }\rightarrow
W,$ such that the induced $R$-bilinear map
\begin{equation*}
V\times W\rightarrow R,\text{ }(v,w)\mapsto \text{ }<v,w>:=\kappa
_{P}(v)(w)=\chi _{P}(w)(v)
\end{equation*}
has the property
\begin{equation*}
<\xi (v),w^{\prime }>=<v,\theta (w^{\prime })>\text{ for all }v\in V\text{
and }w^{\prime }\in W^{\prime }.
\end{equation*}
We say an $R$-pairing $P=(V,W)$ satisfies the $\alpha $\emph{-condition}%
\textbf{\ }(or $P$ is an $\alpha $\emph{-pairing}), if for every $R$-module $%
M$ the following mapping is injective:
\begin{equation}
\alpha _{M}^{P}:M\otimes _{R}W\rightarrow \mathrm{Hom}_{R}(V,M),\text{ }\sum
m_{i}\otimes w_{i}\mapsto \lbrack v\mapsto \sum m_{i}<v,w_{i}>].  \label{alp}
\end{equation}

We say an $R$-module $W$ satisfies the $\alpha $\emph{-condition}, if $%
(W^{\ast },W)$ satisfies the $\alpha $-condition (equivalently, if $_{R}W$
is \emph{locally projective} in the sense of Zimmermann-Huisgen
\cite[Theorem 2.1]{Z-H76}, \cite[Theorem 3.2]{Gar76}). With $\mathcal{P}$ we
denote the category of $R$-pairing with morphisms of pairings described
above and with $\mathcal{P}^{\alpha }\subseteq \mathcal{P}$ the full
subcategory of $R$-pairings satisfying the $\alpha $-condition.
\end{punto}

\begin{remark}
\label{flat}(\cite[Remark 2.2]{Abu1}) Let $P=(V,W)$ be an $\alpha $-pairing.
Then $_{R}W$ is flat and $R$-cogenerated. If moreover $_{R}W$ is finitely
presented or $R$ is perfect, then $_{R}W$ turns to be projective.
\end{remark}

\begin{lemma}
\label{q-2}\emph{(\cite[Lemma 2.3]{Abu1})} Let $P=(V,W)\in \mathcal{P}%
^{\alpha }.$ For every $R$-module $M$ and every $R$-submodule $N\subset M$
we have for $\sum m_{i}\otimes w_{i}\in M\otimes _{R}W:$%
\begin{equation}
\sum m_{i}\otimes w_{i}\in N\otimes _{R}W\Leftrightarrow \sum m_{i}<v,w_{i}>%
\text{ }\in N\text{ for all }v\in V.  \label{NotD}
\end{equation}
\end{lemma}

\begin{punto}
\textbf{Measuring }$R$-\textbf{pairings.}\label{mP}\textbf{\ }For an $R$%
-coalgebra $C$ and an $R$-algebra $A$ (not necessarily with unity) we call
an $R$-pairing $P=(A,C)$ a \emph{measuring }$R$\emph{-pairing}, if the
induced mapping $\kappa _{P}:A\rightarrow C^{\ast }$ is an $R$-algebra
morphism. In this case $C$ is an $A$-bimodule through the left and the right
$A$-actions
\begin{equation}
a\rightharpoonup c:=\sum c_{1}<a,c_{2}>\text{ and }c\leftharpoonup a:=\sum
<a,c_{1}>c_{2}\text{ for all }a\in A,\text{ }c\in C.  \label{C-r}
\end{equation}
Let $(A,C)$ and $(B,D)$ be measuring $R$-pairings ($A$ and $B$ not
necessarily with unities). Then we say an $R$-pairings morphism
$(\xi
,\theta ):(B,D)\rightarrow (A,C)$ is a \emph{morphism of measuring }$R$\emph{%
-pairings,} if $\xi :A\rightarrow B$ is an $R$-algebra morphism and $\theta
:D\rightarrow C$ is an $R$-coalgebra morphism. The category of measuring $R$%
-pairings and morphisms described above will be denoted by $\mathcal{P}_{m}.$
If $P=(A,C)$ is a measuring $R$-pairing, $D\subset C$ is a (pure) $R$%
-subcoalgebra and $I\vartriangleleft A$ is an ideal with $<I,D>=0,$ then $%
Q:=(A/I,D)$ is a measuring $R$-pairing, $(\pi _{I},\iota
_{D}):(A/I,D)\rightarrow (A,C)$ is a morphism in $\mathcal{P}_{m}$ and we
call $Q\subset P$ a (pure) \emph{measuring }$R$\emph{-subpairing}. With $%
\mathcal{P}_{m}^{\alpha }\subset \mathcal{P}_{m}$ we denote the \emph{full }%
subcategory of measuring $R$-pairings satisfying the $\alpha $-condition.
Obviously $\mathcal{P}_{m}^{\alpha }\subset \mathcal{P}_{m}$ is closed under
pure measuring $R$-subpairings.
\end{punto}

\section*{Rational modules}

\begin{punto}
\label{rat-dar}Let $P=(A,C)$ be a measuring $\alpha $-pairing ($A$ not
necessarily with unity). Let $M$ be an $A$-faithful left $A$-module and
consider the injective canonical $A$-linear mapping $\rho _{M}:M\rightarrow
\mathrm{Hom}_{R}(A,M).$ We put $\mathrm{Rat}^{C}(_{A}M):=\rho
_{M}^{-1}(M\otimes _{R}C).$ If $\mathrm{Rat}^{C}(_{A}M)=M,$ then $M$ is said
to be $C$\emph{-rational} and we define
\begin{equation*}
\varrho _{M}:=(\alpha _{M}^{P})^{-1}\circ \rho _{M}:M\rightarrow M\otimes
_{R}C.
\end{equation*}
For $m\in \mathrm{Rat}^{C}(_{A}M)$ with $\varrho
_{M}(m)=\sum\limits_{i=1}^{k}m_{i}\otimes c_{i}$ we call $\left\{
(m_{i},c_{i})\right\} _{i=1}^{k}\subset M\times C$ a \emph{rational system}
for $m.$ With $\mathrm{Rat}^{C}(_{A}\widetilde{\mathcal{M}})\subseteq $ $_{A}%
\widetilde{\mathcal{M}}$ we denote the \emph{full} subcategory of $C$%
-rational left $A$-modules. The \emph{full }subcategory of $C$-rational
right $A$-modules $^{C}\mathrm{Rat}(\widetilde{\mathcal{M}}_{A})\subseteq
\widetilde{\mathcal{M}}_{A}$ is analogously defined (we will show in Theorem
\ref{equal} that every $C$-rational left, respectively right, $A$-module is
unital).
\end{punto}

\qquad As a preparation for the proof of the main results in this section
(Theorems \ref{equal} and \ref{cor-dicht}) and to make the paper more
self-contained we begin with some technical lemmas.

\begin{lemma}
\label{clos}Let $P=(A,C)$ be a measuring $\alpha $-pairing \emph{(}A not
necessarily with unity\emph{)}. For every $A$-faithful left $A$-module $M$
we have:

\begin{enumerate}
\item  $\mathrm{Rat}^{C}(_{A}M)\subseteq M$ is an $A$-submodule.

\item  For every $A$-submodule $N\subset M$ we have $\mathrm{Rat}%
^{C}(_{A}N)=N\cap \mathrm{Rat}^{C}(_{A}M).$

\item  $\mathrm{Rat}^{C}(\mathrm{Rat}^{C}(_{A}M))=\mathrm{Rat}^{C}(_{A}M).$

\item  For every $L\in $ $_{A}\widetilde{\mathcal{M}}$ and $f\in $ $\mathrm{%
Hom}_{A-}(M,L)$ we have $f(\mathrm{Rat}^{C}(_{A}M))\subseteq \mathrm{Rat}%
^{C}(_{A}L).$
\end{enumerate}
\end{lemma}

\begin{Beweis}
\begin{enumerate}
\item  Let $b\in A$ and $m\in \mathrm{Rat}^{C}(_{A}M)$ with rational system $%
\left\{ (m_{i},c_{i})\right\} _{i=1}^{k}\subset M\times C.$ Then we have for
arbitrary $a\in A:$
\begin{equation*}
a(bm)=(ab)m=\sum\limits_{i=1}^{k}m_{i}<ab,c_{i}>=\sum%
\limits_{i=1}^{k}m_{i}<a,bc_{i}>
\end{equation*}
and so $bm\in \mathrm{Rat}^{C}(_{A}M)$ with rational system $\left\{
(m_{i},bc_{i})\right\} _{i=1}^{k}\subset M\times C.$

\item  Clearly $\mathrm{Rat}^{C}(_{A}N)\subseteq N\cap \mathrm{Rat}%
^{C}(_{A}M).$ On the other hand take $n\in N\cap \mathrm{Rat}^{C}(_{A}M)$
with rational system $\left\{ (m_{i},c_{i})\right\} _{i=1}^{k}\subset
M\times C.$ Then for arbitrary $a\in A$ we have $\sum%
\limits_{i=1}^{k}m_{i}<a,c_{i}>=an\in N$ and so $n\in \mathrm{Rat}^{C}(_{A}N)
$ by Lemma \ref{q-2}.

\item  Follows from 1. and 2.

\item  Let $f:M\rightarrow L$ be a morphism of $A$-faithful left $A$-modules
and take $m\in \mathrm{Rat}^{C}(_{A}M)$ with rational system $\left\{
(m_{i},c_{i})\right\} _{i=1}^{k}\subset M\times C.$ Then for arbitrary $a\in
A$ we have
\begin{equation*}
af(m)=f(am)=f(\sum\limits_{i=1}^{k}m_{i}<a,c_{i}>)=\sum%
\limits_{i=1}^{k}f(m_{i})<a,c_{i}>,
\end{equation*}
i.e. $f(m)\in \mathrm{Rat}^{C}(_{A}L)$ with rational system $\left\{
(f(m_{i}),c_{i})\right\} _{i=1}^{k}\subset L\times C.\blacksquare $
\end{enumerate}
\end{Beweis}

\begin{lemma}
\label{co-rat}Let $P=(A,C)\in \mathcal{P}_{m}$ \emph{(}$A$ not necessarily
with unity\emph{)}.

\begin{enumerate}
\item  If $(M,\varrho _{M})$ is a right $C$-comodule, then $M$ is a left $A$%
-module through
\begin{equation}
\xymatrix{ \rho _{M} : M \ar[rr]^(.45){\alpha_M ^{P} \circ \varrho_M} & &
{\rm Hom} _R (A,M) }  \label{mod-st}
\end{equation}
If $M$ is counital and $A$ has unity, then $_{A}M$ is unital \emph{(}and $A$%
-faithful\emph{)}.

\item  Let $(M,\varrho _{M}),(N,\varrho _{N})$ be right $C$-comodules and
consider the induced left $A$-module structures $(M,\rho _{M}),(N,\rho _{N})$
as in \emph{(\ref{mod-st})}. If $f:M\rightarrow N$ is $C$-colinear, then $f$
is $A$-linear.

\item  Let $N$ be a right $C$-comodule, $K\subset N$ be a right $C$%
-subcomodule and consider the induced left $A$-module structures $(N,\rho
_{N}),$ $(K,\rho _{K})$ as in \emph{(\ref{mod-st})}. Then $K\subset N$ is an
$A$-submodule.
\end{enumerate}
\end{lemma}

\begin{Beweis}
\begin{enumerate}
\item  Set $P\otimes P:=(A\otimes _{R}A,C\otimes _{R}C)$ and consider the
following diagram with commutative trapezoids (where $\zeta ^{l}$ the
isomorphism given by $\zeta ^{l}(\delta )(a\otimes b):=\delta (b)(a)$):
\begin{equation}
\xymatrix{M \ar[rrr]^(.45){\rho_M} \ar[ddd]_{\rho_M} \ar@{=}[dr] & & & {\rm
Hom}_R (A,M) \ar[ddd]^{(\mu,M)} \\ { } & M \ar[d]_(.45){\varrho_M}
\ar[r]^(.45){\varrho_M} & M \otimes_R C \ar[ur]_(.45){\alpha_M ^{P}}
\ar[d]^(.45){id_M \otimes_R \Delta} & \\ & M \otimes_R C
\ar[dl]_(.45){\alpha_M ^{P}} \ar[r]_(.45){\varrho_M \otimes_R id_C} &
\ar[dr]^(.45){\alpha _{M}^{P \otimes P}} M \otimes_R C \otimes_R C & \\ {\rm
Hom}_R (A,M) \ar[rr]_(.45){(A,\rho _M)} & & {\rm Hom}_R (A,{\rm Hom}_R
(A,M)) \ar[r]_(.55){{\zeta}^l} & {\rm Hom}_R (A \otimes_R A,M) }
\label{comod}
\end{equation}
By assumption the internal rectangle is commutative and consequently the
outer rectangle is commutative, i.e. $(M,\rho _{M})$ is a left $A$-module.

If $M$ is counital and $A$ has unity, then for every $m\in M:$
\begin{equation*}
1_{A}m=\varepsilon _{C}m=\sum m_{<0>}\varepsilon _{C}(m_{<1>})=m,
\end{equation*}
i.e. $_{A}M$ is unital \emph{(}and $A$-faithful\emph{)}.

\item  Consider the diagram
\begin{equation}
\xymatrix{M \ar[rrrr]^(.45){f} \ar[dd]_(.45){\varrho_M}
\ar[dr]_(.45){\rho_M} & & & & N \ar[dd]^(.45){\varrho_N}
\ar[dl]^(.45){\rho_N} \\ { } & {\rm Hom} _R(A,M) \ar[rr]^(.45){(A,f)} & &
{\rm Hom} _R(A,N) & \\ M \otimes_R C \ar[ur]^(.45){\alpha_M ^{P}}
\ar[rrrr]_(.45){f \otimes id_C} & & & & \ar[ul]_(.45){\alpha_N ^P} N
\otimes_R C & }  \label{co-lin}
\end{equation}
The lower trapezoid is obviously commutative. Moreover both triangles are
commutative by the definition of $\rho _{M}$ and $\rho _{N}$ (\ref{mod-st}).
If $f$ is $C$-colinear, then the outer rectangle is commutative and
consequently the upper trapezoid is commutative, i.e. $f$ is $A$-linear.

\item  Trivial.$\blacksquare $
\end{enumerate}
\end{Beweis}

\begin{lemma}
\label{rat-co}Let $P=(A,C)\in \mathcal{P}_{m}^{\alpha }$ \emph{(}$A$ not
necessarily with unity\emph{)}.

\begin{enumerate}
\item  If $(M,\rho _{M})\in $ $_{A}\widetilde{\mathcal{M}}$ is $C$-rational,
then $M$ is a \emph{counital }right $C$-comodule through
\begin{equation}
\xymatrix{ \varrho _{M} : M \ar[rr]^{(\alpha_M ^{P}) ^{-1} \circ \rho_M} & &
M \otimes_R C}  \label{comod-st}
\end{equation}

\item  Let $(M,\rho _{M}),(N,\rho _{N})\in $ $_{A}\widetilde{\mathcal{M}}$
be $C$-rational an consider the induced right $C$-comodule structures $%
(M,\varrho _{M}),$ $(N,\varrho _{N})$ as in \emph{(\ref{comod-st}).} Then $%
\mathrm{Hom}^{C}(M,N)=\mathrm{Hom}_{A-}(M,N).$

\item  Let $(N,\rho _{N})\in $ $_{A}\widetilde{\mathcal{M}}$ be $C$-rational
and consider the{\normalsize \ induced }right $C$-comodule structure $%
(N,\varrho _{N})$ as in \emph{(\ref{comod-st}).} If $K\subset N$ is an $A$%
-submodule, then $K$ is a counital right $C$-subcomodule and moreover $%
\varrho _{K}=(\varrho _{N})_{|_{K}}.$
\end{enumerate}
\end{lemma}

\begin{Beweis}
\begin{enumerate}
\item  If $(M,\rho _{M})$ is $C$-rational, then $\rho _{M}(M)\subset \alpha
_{M}^{P}(M\otimes _{R}C)$ (by definition). Moreover $\alpha _{M}^{P}$ is
injective, hence $\varrho _{M}:=(\alpha _{M}^{P})^{-1}\circ \rho
_{M}:M\rightarrow M\otimes _{R}C$ is well defined and we have the
commutative diagram
\begin{equation*}
\xymatrix{ {\rm Hom} _R(A,M) & \\ M \ar@{^{(}->}[u]^{\rho_M}
\ar@{.>}[r]_(.45){\varrho_M} & M \otimes C \ar@{^{(}->}[ul]_{\alpha_M ^{P}} }
\end{equation*}
The right trapezoid in diagram (\ref{comod}) is obviously commutative and by
definition of $\varrho _{M}$ (\ref{comod-st}) all other trapezoids are
commutative. By assumption $M$ is a left $A$-module and so the outer
rectangle is also commutative. By \cite[Lemma 2.8]{Abu1} $\alpha
_{M}^{P\otimes P}$ is injective and consequently the internal rectangle is
commutative, i.e. $(M,\varrho _{M})$ is a right $C$-comodule. Moreover, $%
\rho _{M}$ and $\alpha _{M}^{P}$ are by assumption injective and so $\varrho
_{M}:=\alpha _{M}^{P}\circ \rho _{M}$ is injective, i.e. $M$ is counital.

\item  Let $M,N\in \mathrm{Rat}^{C}(_{A}\mathcal{M})$ and $f:M\rightarrow N$
be $A$-linear. The lower trapezoid in diagram (\ref{co-lin}) is obviously
commutative and by definition of $\varrho _{M},\varrho _{N}$ all triangles
are commutative. If $f$ is $A$-linear then, by the injectivity of $\alpha
_{N}^{P},$ the upper trapezoid is commutative and consequently the outer
triangle is commutative, i.e. $f$ is $C$-colinear. So $\mathrm{Hom}%
_{A-}(M,N)\subseteq \mathrm{Hom}^{C}(M,N)$ and the equality follows from
Lemma \ref{co-rat} (2).

\item  Let $(N,\rho _{N})$ be a $C$-rational left $A$-module. If $K\subset N$
is an $A$-submodule, then by Lemma \ref{clos} (2) $\mathrm{Rat}^{C}(K)=K\cap
\mathrm{Rat}^{C}(M)=K,$ i.e. $K$ is a $C$-rational left $A$-module. By (1)
it follows that $K$ is a counital right $C$-comodule through some $R$-linear
map $\varrho _{K}:K\rightarrow K\otimes _{R}C.$ Moreover $K\overset{\iota
_{K}}{\hookrightarrow }N$ is by assumption $A$-linear and so $C$-colinear
(by 2.), i.e. $K\subset N$ is a $C$-subcomodule. By remark \ref{flat} $_{R}C$
is flat and so $\varrho _{K}=(\varrho _{N})_{|_{K}}.\blacksquare $
\end{enumerate}
\end{Beweis}

\begin{remark}
\label{C-subgenerator}Let $C$ be an $R$-coalgebra and $(N,\varrho _{N})$ be
an arbitrary right $C$-comodule. Let $R^{(\Lambda )}\overset{\pi }{%
\longrightarrow }N\longrightarrow 0$ be a free representation of $N$ in $%
\mathcal{M}_{R}.$ Then
\begin{equation*}
C^{(\Lambda )}\simeq R^{(\Lambda )}\otimes _{R}C\overset{\pi \otimes id}{%
\longrightarrow }N\otimes _{R}C\longrightarrow 0
\end{equation*}
is an epimorphism in $\mathcal{M}^{C}.$ Moreover the injective comodule
structure map $\varrho _{N}:N\rightarrow N\otimes _{R}C$ is $C$-colinear,
i.e. $N$ is a $C$-subcomodule of the $C$-generated $C$-comodule $N\otimes
_{R}C$ and so $C$-subgenerated. Since $N\in \mathcal{M}^{C}$ is arbitrary,
we conclude that $C$ is a subgenerator in $\mathcal{M}^{C}.$
\end{remark}

\begin{punto}
\label{End(C)}For every $R$-coalgebra $C$ we have an $R$-algebra
isomorphism:
\begin{equation*}
\Psi :(C^{\ast },\star )\rightarrow (\mathrm{End}^{C}(C,C)^{op},\circ ),%
\text{ }f\mapsto \lbrack c\mapsto \sum f(c_{1})c_{2}]
\end{equation*}
with inverse $\Phi :g\mapsto \varepsilon \circ g.$ Analogously $(^{C}\mathrm{%
End}(C,C),\circ )\simeq (C^{\ast },\star )$ as $R$-algebras. If $P=(A,C)\in
\mathcal{P}_{m}^{\alpha },$ then we have isomorphisms of $R$-algebras:
\begin{equation*}
C^{\ast }\simeq \text{ }^{C}\mathrm{End}(C)=\mathrm{End}(C_{C^{\ast }})=%
\mathrm{End}(C_{\mathrm{End}^{C}(C)^{op}})=\mathrm{End}(C_{\mathrm{End}%
(_{A}C)^{op}}):=\mathrm{Biend}(_{A}C)
\end{equation*}
and
\begin{equation*}
C^{\ast }\simeq \mathrm{End}^{C}(C)^{op}=\mathrm{End}(_{C^{\ast }}C)^{op}=%
\mathrm{End}(_{^{C}\mathrm{End}(C)}C)^{op}=\mathrm{End}(_{\mathrm{End}%
(C_{A})}C)^{op}:=\mathrm{Biend}(C_{A}),
\end{equation*}
where $\mathrm{Biend}(_{A}C)\ $and $\mathrm{Biend}(C_{A})$ are the \emph{%
biendomorphism rings }of $_{A}C$ and $C_{A},$ respectively (compare
\cite[6.4]{Wis88}).
\end{punto}

\qquad We are now ready to prove the main result in this section:

\begin{theorem}
\label{equal}Let $P=(A,C)\in \mathcal{P}_{m},$ $A$ not necessarily with
unity. If $_{R}C$ is locally projective and $\kappa _{P}(A)\subset C^{\ast }$
is dense with respect to the finite topology on $C^{\ast }\hookrightarrow
C^{C}$, then every right \emph{(}respectively left\emph{)} $C$-comodule is a
\emph{unital} left \emph{(}respectively right\emph{)} $A$-module and we have
category isomorphisms
\begin{equation}
\begin{tabular}{lllllll}
$\mathcal{M}^{{C}}$ & $\simeq $ & $\mathrm{Rat}^{{C}}(_{A}\widetilde{%
\mathcal{M}})$ & $=$ & $\mathrm{Rat}^{{C}}(_{A}\mathcal{M})$ & $=$ & $\sigma
\lbrack _{{A}}${$C$}$]$ \\
& $\simeq $ & $\mathrm{Rat}^{{C}}(_{C^{\ast }}\widetilde{\mathcal{M}})$ & $=$
& $\mathrm{Rat}^{{C}}(_{C^{\ast }}\mathcal{M})$ & $=$ & $\sigma \lbrack
_{C^{\ast }}${$C$}$]$%
\end{tabular}
\label{iso-r}
\end{equation}
and
\begin{equation}
\begin{tabular}{lllllll}
$^{{C}}\mathcal{M}$ & $\simeq $ & $^{{C}}\mathrm{Rat}(\widetilde{\mathcal{M}}%
_{A})$ & $=$ & $^{{C}}\mathrm{Rat}(\mathcal{M}_{A})$ & $=$ & $\sigma \lbrack
${$C$}$_{{A}}]$ \\
& $\simeq $ & $^{{C}}\mathrm{Rat}(\widetilde{\mathcal{M}}_{C^{\ast }})$ & $=$
& $^{{C}}\mathrm{Rat}(\mathcal{M}_{C^{\ast }})$ & $=$ & $\sigma \lbrack ${$C$%
}$_{C^{\ast }}].$%
\end{tabular}
\label{iso-l}
\end{equation}
\end{theorem}

\begin{Beweis}
We prove the category isomorphisms (\ref{iso-r}). The isomorphisms of
categories (\ref{iso-l}) follow by symmetry.

\textbf{Step 1. }$\mathcal{M}^{C}\simeq \mathrm{Rat}^{C}(_{A}\widetilde{%
\mathcal{M}}).$

Since $_{R}C$ satisfies the $\alpha $-condition and $\kappa _{P}(A)\subseteq
C^{\ast }$ is dense, it follows by \cite[Proposition 2.4 (2)]{Abu1} that $%
P\in \mathcal{P}_{m}^{\alpha }.$ For every counital $(M,\varrho _{M})\in
\mathcal{M}^{C}$ we conclude that $\rho _{M}:=\alpha _{M}^{P}\circ \varrho
_{M}$ is injective, i.e. the induced left $A$-module is $A$-faithful. By
Lemmas \ref{co-rat} and \ref{rat-co} we have covariant functors
\begin{equation}
\begin{tabular}{llllllll}
$_{A}(-):$ & $\mathcal{M}^{C}$ & $\rightarrow $ & $\mathrm{Rat}^{C}(_{A}%
\widetilde{\mathcal{M}}),$ & $(-)^{C}:$ & $\mathrm{Rat}^{C}(_{A}\widetilde{%
\mathcal{M}})$ & $\rightarrow $ & $\mathcal{M}^{C},$ \\
& $(M,\varrho _{M})$ & $\mapsto $ & $(M,\alpha _{M}^{P}\circ \varrho _{M}),$
&  & $(M,\rho _{M})$ & $\mapsto $ & $(M,(\alpha _{M}^{P})^{-1}\circ \rho
_{M}),$%
\end{tabular}
\label{cov-fund}
\end{equation}
acting as the identity on morphisms. Obviously
\begin{equation*}
(-)^{C}\circ \text{ }_{A}(-)\simeq id_{\mathcal{M}{\normalsize ^{C}}}\text{
and }_{A}(-)\circ (-)^{C}\simeq id_{\mathrm{Rat}^{C}(_{A}\widetilde{\mathcal{%
M}})},
\end{equation*}
i.e. $\mathcal{M}^{C}\simeq \mathrm{Rat}^{C}(_{A}\widetilde{\mathcal{M}}).$

\textbf{Step 2. }$\mathrm{Rat}^{C}(_{A}\widetilde{\mathcal{M}})=\mathrm{Rat}%
^{C}(_{A}\mathcal{M}).$

Let $(N,\rho _{N})\in \mathrm{Rat}^{C}(_{A}\widetilde{\mathcal{M}})$ and $%
n\in N$ with $\varrho _{N}(n)=\sum\limits_{i=1}^{k}n_{i}\otimes c_{i}.$ By
assumption $\kappa _{P}(A)\subset C^{\ast }$ is dense and so there exists
some $a\in A,$ so that $\kappa _{P}(a)(c_{i})=\varepsilon (c_{i})$ for $%
i=1,...,k.$ Hence
\begin{equation*}
n=\sum\limits_{i=1}^{k}n_{i}\varepsilon
(c_{i})=\sum\limits_{i=1}^{k}n_{i}<a,c_{i}>=an\in N\text{ (i.e. }_{A}N\text{
is \emph{unital}).}
\end{equation*}

\textbf{Step 3.} $\mathcal{M}^{C}=\sigma \lbrack _{{A}}{C}].$

By Remark \ref{flat} $_{R}C$ is flat and so $\mathcal{M}^{C}$ is a
Grothendieck category (e.g. \cite[3.13]{BW03}). Moreover by the previous
lemmas and Remark \ref{C-subgenerator} $\mathcal{M}^{C}\subseteq \sigma
\lbrack _{A}C]$ is a full subcategory of $_{A}\mathcal{M}.$ The equality
follows now by the fact that $\sigma \lbrack _{A}C]\subseteq $ $_{A}\mathcal{%
M}$ is the smallest Grothendieck full subcategory of $_{A}\mathcal{M}$ that
contains $C.$

\textbf{Step 4}{\normalsize . }$C^{\ast }$ has unity $\varepsilon _{C}$ and
by assumption $(C^{\ast },C)\in \mathcal{P}_{m}^{\alpha },$ so the proof
above can be repeated to get
\begin{equation*}
\mathcal{M}{\normalsize ^{C}}\simeq \mathrm{Rat}^{C}(_{C^{\ast }}\widetilde{%
\mathcal{M}})=\mathrm{Rat}^{C}(_{C^{\ast }}\mathcal{M})=\sigma \lbrack
_{C^{\ast }}C].\blacksquare
\end{equation*}
\end{Beweis}

\begin{theorem}
\label{cor-dicht}For a measuring $R$-pairing $P=(${$A$}$,${$C$}$)$\ \emph{(}$%
A$ not necessarily with unity\emph{)}, the following are equivalent:

\begin{enumerate}
\item  $P$ satisfies the $\alpha $-condition;

\item  $\mathcal{M}^{{C}}\simeq \sigma \lbrack _{{A}}${$C$}$]=\sigma \lbrack
_{C^{\ast }}${$C$}$];$

\item  $_{R}C$ is locally projective and $\kappa _{P}(A)\subseteq C^{\ast }$
is dense;

\item  $^{{C}}\mathcal{M}\simeq \sigma \lbrack ${$C$}$_{{A}}]=\sigma \lbrack
${$C$}$_{C^{\ast }}].$
\end{enumerate}
\end{theorem}

\begin{Beweis}
1. $\Rightarrow $ 2. Follows from Theorem \ref{equal}.

2. $\Rightarrow $ 3. By assumption we have
\begin{equation*}
\sigma \lbrack _{{A}}{C}]=\sigma \lbrack _{{C}^{\ast }}{C}]=\sigma \lbrack _{%
\mathrm{Biend}(_{A}C)}{C}]
\end{equation*}
and the density of $\kappa _{P}(A)\subseteq C^{\ast }$ follows by Lemma \ref
{dense}. The proof of ``$\mathcal{M}^{C}=\sigma \lbrack _{C^{\ast }}C]$ $%
\Rightarrow $ $_{R}C$ locally projective'' follows from \cite[3.5]{Wis02}
(which appeared also as \cite[4.3]{BW03}).

3. $\Rightarrow $ 1. Follows from general theory of dual pairings over rings
(e.g. \cite[Proposition 2.4 (2)]{Abu1}).

1. $\Leftrightarrow $ 4. Follows by symmetry.$\blacksquare $
\end{Beweis}

\begin{example}
An interesting example for a measuring pairing for which Theorem \ref
{cor-dicht} applies is $P:=(C^{\Box },C),$ where $C$ is a locally projective
$R$-coalgebra and $C^{\Box }:=\mathrm{Rat}^{C}(_{C^{\ast }}C).$
\end{example}

\section*{Dual coalgebras}

\begin{definition}
Let $A$ be an $R$-algebra and consider the class of $R$-cofinite $A$-ideals $%
\mathcal{K}_{A}.$ For every class $\mathcal{F}$ of $R$-cofinite $A$-ideals
we define the \emph{set}
\begin{equation}
A_{\mathcal{F}}^{\circ }:=\{f\in A^{\ast }\mid f(I)=0\text{ for some }I\in
\mathcal{F}\}.  \label{A0-F}
\end{equation}

\begin{enumerate}
\item  A filter $\frak{F}=\{I_{\lambda }\}_{\Lambda }$ consisting of $R$%
-cofinite $A$-ideals will be called

an $\alpha $\emph{-filter}, if the $R$-pairing $(A,A_{\frak{F}}^{\circ })$
satisfies the $\alpha $-condition;

\emph{cofinitary}, if for every $I_{\lambda }\in \frak{F}$\ there exists $%
I_{\varkappa }\subset I_{\lambda }$ for some $\varkappa \in \Lambda ,$ such
that $A/I_{\varkappa }$ is finitely generated and projective in $\mathcal{M}%
_{R};$

\emph{cofinitely }$R$\emph{-cogenerated,} if $A/I$ is $R$-cogenerated for
every $I\in \frak{F}.$

\item  We call $A:$

an $\alpha $\emph{-algebra}, if $\mathcal{K}_{A}$ is an $\alpha $-filter;

\emph{cofinitary}, if $\mathcal{K}_{A}$ is a cofinitary filter;

\emph{cofinitely }$R$\emph{-cogenerated}, if $A/I$ is $R$-cogenerated for
every $I\in \mathcal{K}_{A}.$
\end{enumerate}
\end{definition}

\begin{Notation}
With $\mathbf{Cog}_{R}$ (respectively $\mathbf{Big}_{R},$ $\mathbf{Hopf}_{R}$%
) denote the category of $R$-coalgebras (respectively $R$-bialgebras, Hopf $R
$-algebras) and with $\mathbf{CAlg}_{R}$ (respectively $\mathbf{CCog}_{R}$)
the category of \emph{commutative} $R$-algebras (respectively \emph{%
cocommutative }$R$-coalgebras). With $\mathbf{CBig}_{R}$ (respectively $%
\mathbf{CCBig}_{R}$) we denote the category of commutative (respectively
cocommutative) $R$-bialgebras and with $\mathbf{CHopf}_{R}$ (respectively $%
\mathbf{CCHopf}_{R}$) the category of commutative (respectively
cocommutative) Hopf $R$-algebras.

For two $R$-coalgebras $C,D$ we denote with $\mathrm{Cog}_{R}(C,D)$ the set
of all $R$-coalgebra morphisms from $C$ to $D.$ For two $R$-algebras
(respectively $R$-bialgebras, Hopf $R$-algebras) $H,$ $K$ we denote with $%
\mathrm{Alg}_{R}(H,K)$ (respectively $\mathrm{Big}_{R}(H,K),$ $\mathrm{Hopf}%
_{R}(H,K)$) the set of all $R$-algebra morphisms (respectively $R$-bialgebra
morphisms, Hopf $R$-morphisms) from $H$ to $K.$
\end{Notation}

\begin{remark}
We make the convention that an $R$-bialgebra (respectively a Hopf $R$%
-algebra) is an $\alpha $\emph{-bialgebra }(respectively a \emph{Hopf }$%
\alpha $\emph{-algebra}), is \emph{cofinitary} or is \emph{cofinitely }$R$%
\emph{-cogenerated}, if it is so as an $R$-algebra. With $\mathbf{Big}%
_{R}^{\alpha }\subset \mathbf{Big}_{R}$ (respectively $\mathbf{Hopf}%
_{R}^{\alpha }\subset \mathbf{Hopf}_{R}$) we denote the full subcategory of $%
\alpha $-bialgebras (respectively Hopf $\alpha $-algebras).
\end{remark}

\begin{lemma}
\label{A0=}\emph{(\cite[Proposition 2.6]{AG-TW2000})} Let $R$ be Noetherian
and $A$ be an $R$-algebra. Then
\begin{equation*}
\begin{tabular}{lll}
$A^{\circ }$ & $:=$ & $\{f\in A^{\ast }\mid f(I)=0$ for some $R$-cofinite
ideal $I\vartriangleleft A\};$ \\
& $=$ & $\{f\in A^{\ast }|\text{ }f(I)=0$ for some $R$-cofinite left \emph{(}%
right\emph{)} $A$-ideal$\};$ \\
& $=$ & $\{f\in A^{\ast }|\text{ }Af$ \emph{(}$fA$\emph{)} is f.g. in $%
\mathcal{M}_{R}\}$ \\
& $=$ & $\{f\in A^{\ast }|\text{ }AfA\text{ is f.g. in }\mathcal{M}_{R}\}.$%
\end{tabular}
\end{equation*}
\end{lemma}

\begin{theorem}
\label{R(G)-co}\emph{(\cite[Theorem 3.3.]{Abu2}) }Let $R$ be Noetherian, $A$
be an $R$-algebra and consider $A^{\circ }\subseteq A^{\ast }$ as an $A$%
-bimodule under the left and the right \emph{regular }$A$\emph{-actions}
\begin{equation}
(af)(\widetilde{a})=f(\widetilde{a}a)\text{ and }(fa)(\widetilde{a})=f(a%
\widetilde{a})\text{ for all }a,\widetilde{a}\in A\text{ and }f\in A^{\ast }.
\label{regular}
\end{equation}
For an $A$\emph{-subbimodule }$C\subseteq A^{\circ }$ and $P:=(A,C)$ the
following are equivalent:

\begin{enumerate}
\item  $_{R}C$ is locally projective and $\kappa _{P}(A)\subset C^{\ast }$
is dense;

\item  $_{R}C$ satisfies the $\alpha $-condition and $\kappa _{P}(A)\subset
C^{\ast }$ is dense;

\item  $(A,C)$ is an $\alpha $-pairing;

\item  $C\subset R^{A}$ is pure;

\item  $C$ is an $R$-coalgebra and $(A,C)\in \mathcal{P}_{\alpha }^{m}.$

If $R$ is a QF Ring, then these are moreover equivalent to

\item  $_{R}C$ is projective.
\end{enumerate}
\end{theorem}

\begin{corollary}
\emph{(\cite[Corollary 3.16]{Abu2})} Let $A$ be an $R$-algebra and $\frak{F}%
\ $be a filter consisting of $R$-cofinite $A$-ideals. If $R$ is Noetherian
and $\frak{F}$ is an $\alpha $-filter, or if $\frak{F}$ is cofinitary then
we have isomorphisms of categories
\begin{equation*}
\begin{tabular}{lllll}
$\mathcal{M}^{A_{\frak{F}}^{\circ }}$ & $\simeq $ & $\mathrm{Rat}^{A_{\frak{F%
}}^{\circ }}(_{A}\mathcal{M})$ & $=$ & $\sigma \lbrack _{A}A_{\frak{F}%
}^{\circ }]$ \\
& $\simeq $ & $\mathrm{Rat}^{A_{\frak{F}}^{\circ }}(_{A_{\frak{F}}^{\circ
\ast }}\mathcal{M})$ & $=$ & $\sigma \lbrack _{A_{\frak{F}}^{\circ \ast }}A_{%
\frak{F}}^{\circ }]$%
\end{tabular}
\text{ }\&\text{ }
\begin{tabular}{lllll}
$^{A_{\frak{F}}^{\circ }}\mathcal{M}$ & $\simeq $ & $^{A_{\frak{F}}^{\circ }}%
\mathrm{Rat}(\mathcal{M}_{A})$ & $=$ & $\sigma \lbrack A_{\frak{F}A}^{\circ
}]$ \\
& $\simeq $ & $^{A_{\frak{F}}^{\circ }}\mathrm{Rat}(\mathcal{M}_{A_{\frak{F}%
}^{\circ \ast }})$ & $=$ & $\sigma \lbrack A_{\frak{F}A_{\frak{F}}^{\circ
\ast }}^{\circ }].$%
\end{tabular}
\end{equation*}
\end{corollary}

\section{The cotensor functor}

\qquad Dual to the \emph{tensor product} of modules, J. Milnor and J. Moore
introduced in \cite{MM65} the \emph{cotensor product} of comodules. For a
closer look on the properties of the cotensor product over arbitrary
(commutative) base rings the interested reader may refer to \cite{Guz85}
(and \cite{Alt99}).

\begin{punto}
Let $C$ be an $R$-coalgebra, $(M,\varrho _{M})\in \mathcal{M}^{C},$ $%
(N,\varrho _{N})\in $ $^{C}\mathcal{M}$ and consider the $R$-linear mapping
\begin{equation*}
\overline{\varrho }_{M,N}:=\varrho _{M}\otimes id_{N}-id_{M}\otimes \varrho
_{N}:M\otimes _{R}N\rightarrow M\otimes _{R}C\otimes _{R}N.
\end{equation*}
The \emph{cotensor product} of $M\;$and $N$ (denoted with $M\square _{C}N$)
is defined through the exactness of the following sequence in $\mathcal{M}%
_{R}:$%
\begin{equation*}
0\rightarrow M\square _{C}N\rightarrow M\otimes _{R}N\overset{\overline{%
\varrho }_{M,N}}{\rightarrow }M\otimes _{R}C\otimes _{R}N.
\end{equation*}
For $M,M^{\prime }\in \mathcal{M}^{C}$ and $N,N^{\prime }\in $ $^{C}\mathcal{%
M},$ the \emph{cotensor product} of $f\in \mathrm{Hom}^{C}(M,M^{\prime })$
and $g\in $ $^{C}\mathrm{Hom}(N,N^{\prime })$ is defined as the $R$-linear
mapping
\begin{equation*}
f\square _{C}g:M\square _{C}N\rightarrow M^{\prime }\square _{C}N^{\prime },
\end{equation*}
that completes the following diagram commutatively
\begin{equation}
\xymatrix{ 0 \ar[r] & M \Box _C N \ar[r] \ar@{-->}[d]^{f \Box_C g} & M
\otimes_R N \ar[r]^(.4){\overline{\varrho}_{M,N}} \ar[d]^(.45){f \otimes g}
& M \otimes_R C \otimes_R N \ar[d]^(.45){f \otimes id_C \otimes g} \\ 0
\ar[r] & M' \Box_C N' \ar[r] & M' \otimes_R N'
\ar[r]^(.4){\overline{\varrho}_{M',N'}} & M' \otimes_R C \otimes_R N'}
\label{f-ct-g}
\end{equation}
In this way we get the \emph{cotensor functor}
\begin{equation*}
M\square _{C}-:\text{ }^{C}\mathcal{M}\rightarrow \mathcal{M}_{R}\text{ }\;%
\text{\emph{(}respectively }-\square _{C}N:\mathcal{M}^{C}\rightarrow
\mathcal{M}_{R}\text{\emph{)}},
\end{equation*}
which is left exact if $_{R}C$ and $M_{R}$ (respectively $_{R}C$ and $_{R}N$%
) are flat.
\end{punto}

\qquad

\begin{definition}
Let $C$ be a \emph{flat} $R$-coalgebra (hence $\mathcal{M}^{C}$ and $^{C}%
\mathcal{M}$ are \emph{abelian} categories). A right (respectively a left) $C
$-comodule $M\;$is called \emph{coflat}, if the functor $M\square _{C}-:$ $%
^{C}\mathcal{M}\rightarrow \mathcal{M}_{R}$ (respectively $-\square _{C}M:%
\mathcal{M}^{C}\rightarrow \mathcal{M}_{R}$) is exact.
\end{definition}

\begin{lemma}
\label{koass2}\emph{(Compare \cite[Page 127]{Sch90}, \cite[10.6]{BW03})} Let
$C$ be an $R$-coalgebra and $M\in \mathcal{M}^{C},$ $N\in $ $^{C}\mathcal{M}.
$ If $W_{R}$ is flat, then there are isomorphisms of $R$-modules
\begin{equation}
W\otimes _{R}(M\square _{C}N)\simeq (W\otimes _{R}M)\square _{C}N\text{ and }%
(M\square _{C}N)\otimes _{R}W\simeq M\square _{C}(N\otimes _{R}W).
\label{k-ot-assoc}
\end{equation}
\end{lemma}

The following result can easily be derived with the help of Lemma \ref
{koass2}:

\begin{corollary}
\label{cot-com}Let $C,D$ be $R$-coalgebras and $(M,\varrho _{M}^{C},\varrho
_{M}^{D})\in $ $^{C}\mathcal{M}^{D}.$

\begin{enumerate}
\item  Assume $C_{R}$ to be flat. For every $N\in $ $^{D}\mathcal{M},$ $%
M\square _{D}N$ is a left $C$-comodule through
\begin{equation*}
\varrho _{M}^{C}\square _{D}id_{N}:M\square _{D}N\rightarrow (C\otimes
_{R}M)\square _{D}N\simeq C\otimes _{R}(M\square _{D}N).
\end{equation*}

\item  Assume $_{R}D$ to be flat. For every $L\in \mathcal{M}^{C},$ $%
L\square _{C}M$ is a right $D$-comodule through
\begin{equation*}
id_{L}\square _{C}\varrho _{M}^{D}:L\square _{C}M\mapsto L\square
_{C}(M\otimes _{R}D)\simeq (L\square _{C}M)\otimes _{R}D.
\end{equation*}
\end{enumerate}
\end{corollary}

\begin{remark}
\label{ko-flat}(\cite[Lemma II.2.5, Folgerung II.2.6]{Alt99}) Let $C$ be a
\emph{flat} $R$-coalgebra. For every $M\in \mathcal{M}^{C},$ the mapping $%
\varrho _{M}:M\rightarrow M\square _{C}C$ is an isomorphism in $\mathcal{M}%
^{C}$ with inverse $\lambda _{M}:m\otimes c\mapsto m\varepsilon (c)$ and
moreover we have
\begin{equation*}
M\otimes _{R}-\simeq M\square _{C}(C\otimes _{R}-):\mathcal{M}%
_{R}\rightarrow \mathcal{M}_{\mathbb{Z}}.
\end{equation*}
If $M$ is coflat in $\mathcal{M}^{C},$ then $M_{R}$ is flat.
\end{remark}

The Associativity of the cotensor products is not valid in general (see \cite
{GP87}). However we have it in special cases, e.g. :

\begin{lemma}
\label{koass-cot}\emph{(\cite[Folgerung II.3.4.]{Alt99})} Let $C,$ $D$ be
\emph{flat} $R$-coalgebras, $N\in $ $^{D}\mathcal{M},$ $M\in $ $^{C}\mathcal{%
M}^{D}$ and $L\in \mathcal{M}^{C}.$ If $L\in \mathcal{M}^{C}$ \emph{(}or $%
N\in $ $^{D}\mathcal{M}$\emph{)} is coflat, then we have an isomorphism of $R
$-modules
\begin{equation}
(L\square _{C}M)\square _{D}N\simeq L\square _{C}(M\square _{D}N).
\label{kot-koass}
\end{equation}
\end{lemma}

\begin{Notation}
For an $R$-algebra $A$ we denote with $A^{e}:=A\otimes _{R}A^{op}$ the \emph{%
enveloping }$R$\emph{-algebra} of $A.$
\end{Notation}

\begin{lemma}
\label{N-funk}Let $A$ be an $R$-algebra, $M,N\in $ $_{A}\mathcal{M}$ and
consider $A,$ $\mathrm{Hom}_{R}(N,M)$ with the canonical left $A^{e}$-module
structures. Then we have a functorial isomorphism
\begin{equation*}
\mathrm{Hom}_{A^{e}-}(A,\mathrm{Hom}_{R}(N,M))\simeq \mathrm{Hom}_{A-}(N,M).
\end{equation*}
\end{lemma}

\begin{Beweis}
The isomorphism is given by
\begin{equation*}
\Phi _{N,M}:\mathrm{Hom}_{A^{e}-}(A,\mathrm{Hom}_{R}(N,M))\rightarrow
\mathrm{Hom}_{A-}(N,M),\text{ }f\mapsto f(1_{A})
\end{equation*}
with inverse $\Psi _{N,M}:g\mapsto \lbrack a\mapsto ag(-)]].$ One can easily
show that $\Phi _{N,M}$ and $\Psi _{N,M}$ are functorial in $M$ and $%
N.\blacksquare $
\end{Beweis}

In the case of a base field, the cotensor functor is equivalent to a
suitable $\mathrm{Hom}$-functor (e.g. \cite[Proposition 3.1]{AW}). Over
arbitrary ground rings we have

\begin{proposition}
\label{cot=Hom(ot)}Let $P=(A,C)\in \mathcal{P}_{m}$, $(M,\varrho _{M})\in
\mathcal{M}^{C},$ $(N,\varrho _{N})\in $ $^{C}\mathcal{M}$ and consider $A,$
$M\otimes _{R}N$ with the \emph{canonical }left $A^{e}$-module structures.

\begin{enumerate}
\item  If $\alpha _{M\otimes _{R}N}^{P}$ is injective, then we have for $%
\sum m_{i}\otimes n_{i}\in M\otimes _{R}N:$
\begin{equation*}
\sum m_{i}\otimes n_{i}\in M\square _{C}N\Leftrightarrow \sum am_{i}\otimes
n_{i}=\sum m_{i}\otimes n_{i}a\text{ for all }a\in A.
\end{equation*}

\item  If $P\in \mathcal{P}_{m}^{\alpha },$ then we have a functorial
isomorphism
\begin{equation*}
M\square _{C}N\simeq \mathrm{Hom}_{A^{e}-}(A,M\otimes _{R}N).
\end{equation*}
\end{enumerate}
\end{proposition}

\begin{Beweis}
\begin{enumerate}
\item  Let $\alpha _{M\otimes _{R}N}^{P}$ be injective and set $\psi
:=\alpha _{M\otimes _{R}N}^{P}\circ \tau _{(23)}.$ Then
\begin{equation*}
\begin{tabular}{llll}
& $\sum m_{i}\otimes n_{i}\in M\square _{C}N$ &  &  \\
$\Leftrightarrow $ & $\sum m_{i<0>}\otimes m_{i<1>}\otimes n_{i}$ & $=$ & $%
\sum m_{i}\otimes n_{i<-1>}\otimes n_{i<0>},$ \\
$\Leftrightarrow $ & $\psi (\sum m_{i<0>}\otimes m_{i<1>}\otimes n_{i})(a)$
& $=$ & $\psi (\sum m_{i}\otimes n_{i<-1>}\otimes n_{i<0>})(a),$ $\forall $ $%
a\in A$ \\
$\Leftrightarrow $ & $\sum m_{i<0>}<a,m_{i<1>}>\otimes n_{i}$ & $=$ & $\sum
m_{i}\otimes <a,n_{i<-1>}>n_{i<0>},$ $\forall $ $a\in A$ \\
$\Leftrightarrow $ & $\sum am_{i}\otimes n_{i}$ & $=$ & $\sum m_{i}\otimes
n_{i}a,$ $\forall $ $a\in A.$%
\end{tabular}
\end{equation*}

\item  The isomorphism is given through
\begin{equation*}
\gamma _{M,N}:M\square _{C}N\rightarrow \mathrm{Hom}_{A^{e}-}(A,M\otimes
_{R}N),\text{ }m\otimes n\mapsto \lbrack a\mapsto am\otimes n\text{ }%
(=m\otimes na)]
\end{equation*}
with inverse $\beta _{M,N}:f\mapsto f(1_{A}).$ It is easy to see that $%
\gamma _{M,N}$ and $\beta _{M,N}$\ are functorial in $M$ and $N.\blacksquare
$
\end{enumerate}
\end{Beweis}

\begin{lemma}
\label{vm-Alg}\emph{(\cite[15.7]{Wis96}, \cite[II, 4.2, Proposition 2]{Bou74}%
)} Let $A$ be an $R$-algebra, $K,K^{\prime }$ be left $A$-modules, $L$ be an
$R$-module and consider the $R$-linear mapping
\begin{equation}
\upsilon :\mathrm{Hom}_{A-}(K,K^{\prime })\otimes _{R}L\rightarrow \mathrm{%
Hom}_{A-}(K,K^{\prime }\otimes _{R}L),\text{ }h\otimes l\mapsto h(-)\otimes
l.  \label{vm}
\end{equation}

\begin{enumerate}
\item  If $_{R}L$ is flat and $_{A}K$ is finitely generated \emph{(}%
respectively finitely presented\emph{)}, then $\upsilon $ is injective \emph{%
(}respectively bijective\emph{)}.

\item  If $_{A}K$ be $K^{\prime }$-projective and $_{A}K$ is finitely
generated, then $\upsilon $ is bijective.

\item  If $_{A}K$ be $K^{\prime }$-projective and $_{R}L$ is finitely
presented, then $\upsilon $ is bijective.

\item  If $_{R}L\;$is projective \emph{(}respectively finitely generated
projective\emph{)}, then $\upsilon $ is injective \emph{(}respectively
bijective\emph{)}.
\end{enumerate}
\end{lemma}

\begin{punto}
\label{fp}(\cite[3.11]{BW03}) Let $_{R}C\;$be a flat $R$-coalgebra. Let $M$
be a left $C$-comodule and consider the $R$-linear mapping
\begin{equation}
\gamma :M^{\ast }\rightarrow \mathrm{Hom}_{R}(M,C),\text{ }f\mapsto \lbrack
m\mapsto \sum f(m_{<-1>})m_{<0>}].  \label{rh_Mst}
\end{equation}
If $_{R}M$ is finitely presented, then $\mathrm{Hom}_{R}(M,C)\simeq M^{\ast
}\otimes _{R}C$ (see Lemma \ref{vm-Alg}) and $M^{\ast }$ is a right $C$%
-comodule through
\begin{equation}
\varrho _{M^{\ast }}:M^{\ast }\overset{\gamma }{\rightarrow }\mathrm{Hom}%
_{R}(M,C)\simeq M^{\ast }\otimes _{R}C.  \label{vr_Mst}
\end{equation}
If $M$ is a right $C$-comodule and $M_{R}$ is finitely presented, then $%
M^{\ast }$ becomes analogously a left $C$-comodule.
\end{punto}

With the help of Lemmas \ref{N-funk} and \ref{vm-Alg}, the following result
can be derived directly from Proposition \ref{cot=Hom(ot)}:

\begin{corollary}
Let $P=(A,C)\in \mathcal{P}_{m}^{\alpha }.$

\begin{enumerate}
\item  Let $M,N\in $ $\mathcal{M}^{C}.$ If $M_{R}$ is flat and $_{R}N$ is
finitely presented, or $N_{R}$ is finitely generated projective, then we
have functorial isomorphisms
\begin{equation*}
\begin{tabular}{lllll}
$M\square _{C}N^{\ast }$ & $\simeq $ & $\mathrm{Hom}_{A^{e}-}(A,M\otimes
_{R}N^{\ast })$ & $\simeq $ & $\mathrm{Hom}_{A^{e}-}(A,\mathrm{Hom}_{R}(N,M))
$ \\
& $\simeq $ & $\mathrm{Hom}_{A-}(N,M)$ & $=$ & $\mathrm{Hom}^{C}(N,M).$%
\end{tabular}
\end{equation*}

\item  Let $M\in \mathcal{M}^{C},$ $N$ be a $C$-bicomodule and consider $N$
with the induced left $A^{e}$-module structure. Then we have isomorphisms of
$R$-modules
\begin{equation*}
M\square _{C}N\simeq \mathrm{Hom}{\normalsize _{A^{e}-}}(A,M\otimes
_{R}N)\simeq \text{ }M\otimes _{R}\mathrm{Hom}{\normalsize _{A^{e}-}}(A,N),
\end{equation*}
if any one of the following conditions is satisfied:

\begin{enumerate}
\item  $M_{R}$ is flat and $_{A^{e}}A$ is finitely presented \emph{(}e.g. $A$
is an affine $R$-algebra \emph{\cite[23.6]{Wis96})};

\item  $_{A^{e}}A$ is $N$-projective and finitely generated;

\item  $_{A^{e}}A$ is $N$-projective and $M_{R}$ is finitely presented;

\item  $M_{R}$ is finitely generated projective.
\end{enumerate}

\item  Let $N\in $ $^{C}\mathcal{M},$ $M$ be a $C$-bicomodule and consider $M
$ with the induced left $A^{e}$-module structure. Then we have an
isomorphism of $R$-modules
\begin{equation*}
M\square _{C}N\simeq \mathrm{Hom}_{A^{e}-}(A,M\otimes _{R}N)\simeq \mathrm{%
Hom}_{A^{e}-}(A,M)\otimes _{R}N,
\end{equation*}
if any one of the following conditions is satisfied:

\begin{enumerate}
\item  $_{R}N$ is flat and $_{A^{e}}A$ is finitely presented \emph{(}e.g. $A$
is an affine $R$-algebra \emph{\cite[23.6]{Wis96})};

\item  $_{A^{e}}A$ is $M$-projective and finitely generated;

\item  $_{A^{e}}A$ is $M$-projective and $_{R}N$ is finitely presented;

\item  $_{R}N$ is finitely generated projective.
\end{enumerate}
\end{enumerate}
\end{corollary}

\subsection*{Injective comodules}

\qquad

For $P=(A,C)\in \mathcal{P}_{m}^{\alpha }$ we get from \cite[16.3]{Wis88}
the following characterizations of the injective objects in $\mathcal{M}%
^{C}\simeq \mathrm{Rat}^{C}(_{A}\mathcal{M})=\sigma \lbrack _{A}C]:$

\begin{lemma}
\label{inj-Com}Let $P=(A,C)\in \mathcal{P}_{m}^{\alpha }.$ For every $U\in
\mathrm{Rat}^{C}(_{A}\mathcal{M})$ the following are equivalent:

\begin{enumerate}
\item  $U$ is injective in $\mathrm{Rat}^{C}(_{A}\mathcal{M});$

\item  $\mathrm{Hom}^{C}(-,U)\simeq \mathrm{Hom}_{A-}(-,U):\mathrm{Rat}%
^{C}(_{A}\mathcal{M})\rightarrow \mathcal{M}_{R}\ $is exact;

\item  $U$ is $C$-injective in $\mathrm{Rat}^{C}(_{A}\mathcal{M});$

\item  $U$ is $K$-injective for every \emph{(}finitely generated, cyclic%
\emph{)} left $A$-submodule $K\subset C;$

\item  every exact sequence $0\rightarrow U\rightarrow L\rightarrow
N\rightarrow 0$ in $\mathrm{Rat}^{C}(_{A}\mathcal{M})$ splits.

\item  every exact sequence $0\rightarrow U\rightarrow L\rightarrow
N\rightarrow 0$ in $\mathrm{Rat}^{C}(_{A}\mathcal{M}),$ in which $N$ is a
factor module of $C$ \emph{(}or $A$\emph{)} splits.
\end{enumerate}
\end{lemma}

The following Lemma plays an important role in the study of injective
objects in the category $\mathrm{Rat}^{C}(_{A}\mathcal{M}),$ where $(A,C)\in
\mathcal{P}_{m}^{\alpha }:$

\begin{lemma}
\label{inj=cof}Let $(A,C)\in \mathcal{P}_{m}^{\alpha }.$ If $R$ is a QF ring
then a $C$-rational left $A$-module $M,$ with $_{R}M$ flat, is injective in $%
\mathrm{Rat}^{C}(_{A}\mathcal{M})$ if and only if $M$ is coflat in $\mathcal{%
M}^{C}.$
\end{lemma}

\begin{Beweis}
By Theorem \ref{cor-dicht} we have the isomorphism of categories
\begin{equation*}
\sigma \lbrack _{A}C]=\mathrm{Rat}{\normalsize ^{C}(_{A}}\mathcal{M})\simeq
\mathcal{M}{\normalsize ^{C}}
\end{equation*}
and we get the result by \cite[10.12]{BW03}.$\blacksquare $
\end{Beweis}

\begin{lemma}
\label{resp-inj}If $P=(A,C)\in \mathcal{P}_{m}^{\alpha },$ then $-\otimes
_{R}C:\mathcal{M}_{R}\rightarrow \mathrm{Rat}^{C}(_{A}\mathcal{M})$ respects
injective objects.
\end{lemma}

\begin{Beweis}
By Theorem \ref{cor-dicht} $\mathcal{M}^{C}\simeq \mathrm{Rat}^{C}(_{A}%
\mathcal{M})=\sigma \lbrack _{A}C],$ i.e. $\mathrm{Rat}^{C}(_{A}\mathcal{M}%
)\subset $ $_{A}\mathcal{M}$ is a closed subcategory. The \emph{exact}
forgetful functor $\digamma :\mathcal{M}^{C}\rightarrow \mathcal{M}_{R}$ is
left adjoint to $-\otimes _{R}C:\mathcal{M}_{R}\rightarrow \mathcal{M}^{C}$
and the result follows then by \cite[45.6]{Wis88}.$\blacksquare $
\end{Beweis}

\begin{proposition}
\label{sep->h.e.}Let $(A,C)\in \mathcal{P}_{m}^{\alpha }$ and $M\in \mathrm{%
Rat}^{C}(_{A}\mathcal{M}).$

\begin{enumerate}
\item  $M$ is an $A$-submodule of an injective $C$-rational left $A$-module.

\item  Every injective object in $\mathrm{Rat}^{C}(_{A}\mathcal{M})$ is $C$%
-generated.

\item  $M$ is injective in $\mathrm{Rat}^{C}(_{A}\mathcal{M})$ if and only
if there exists an injective $R$-module $X$ for which $_{A}M$ is a direct
summand of $X\otimes _{R}C.$

\item  Let $M$ be injective in $\mathcal{M}_{R}.$ Then $M$ is injective in $%
\mathrm{Rat}^{C}(_{A}\mathcal{M})$ if and only if $\varrho _{M}:M\rightarrow
M\otimes _{R}C$ splits in $_{A}\mathcal{M}.$

\item  Let $R$ be Noetherian. Then $M$ is injective in $\mathrm{Rat}^{C}(_{A}%
\mathcal{M})$ if and only if $M^{(\Lambda )}$ is injective in $\mathrm{Rat}%
^{C}(_{A}\mathcal{M})$ for every index set $\Lambda .$ Moreover, direct
limits of injectives in $\mathrm{Rat}^{C}(_{A}\mathcal{M})$ are injective.

\item  Let $A$ be separable \emph{(}i.e. $_{A^{e}}A$ is projective\emph{)}.
Then $M\in \mathcal{M}^{C}$ is coflat if and only if $M_{R}$ is flat.
\end{enumerate}
\end{proposition}

\begin{Beweis}
\begin{enumerate}
\item  Let $M\in \mathrm{Rat}^{C}(_{A}\mathcal{M})$ and denote with $E(M)$
the injective hull of $M$ in $\mathcal{M}_{R}.$ By Lemma \ref{resp-inj} $%
E(M)\otimes _{R}C$ is injective in $\mathrm{Rat}^{C}(_{A}\mathcal{M}).$
Obviously $(\iota _{M}\otimes _{R}id_{C})\circ \varrho _{M}:M\hookrightarrow
E(M)\otimes _{R}C$ is $A$-linear and the result follows.

\item  Let $(M,\varrho _{M})\in \mathrm{Rat}^{C}(_{A}\mathcal{M})$ be
injective. By Lemma \ref{inj-Com}, there exists an epimorphism of left $A$%
-modules $\beta :M\otimes _{R}C\rightarrow M,$ such that $\beta \circ
\varrho _{M}=id_{M}.$ If $R^{(\Lambda )}\overset{\pi }{\rightarrow }%
M\rightarrow 0$ is a free representation of $M$ in $\mathcal{M}_{R},$ then
we get the following exact sequence in $\mathrm{Rat}^{C}(_{A}\mathcal{M}):$
\begin{equation*}
\xymatrix{C^{(\Lambda )} \simeq R^{(\Lambda )} \otimes _{R} C
\ar[rr]^(.65){\beta \circ (\pi \otimes id_{C})} & & M \ar[r] & 0}.
\end{equation*}

\item  Let $X$ be an injective $R$-module, such that $_{A}M$ is a direct
summand of $X\otimes _{R}C.$ By Lemma \ref{resp-inj} $X\otimes _{R}C$ is
injective in $\mathrm{Rat}^{C}(_{A}\mathcal{M})$ and consequently $M$ is
injective in $\mathrm{Rat}^{C}(_{A}\mathcal{M}).$ On the other hand, let $M$
be injective in $\mathrm{Rat}^{C}(_{A}\mathcal{M})$ and denote with $E(M)$
the injective hull of $M$ in $\mathcal{M}_{R}.$ Then we get an exact
sequence in $\mathrm{Rat}^{C}(_{A}\mathcal{M})$%
\begin{equation}
\xymatrix{ 0 \ar[r] & M \ar[rr]^(.4){(\iota \otimes id_C) \circ \varrho _M}
& & E(M) \otimes_R C}.  \label{M_E(M)otC}
\end{equation}
Now (\ref{M_E(M)otC}) splits in $\mathrm{Rat}^{C}(_{A}\mathcal{M})$ by Lemma
\ref{inj-Com} and the result follows.

\item  Follows from Lemmata \ref{inj-Com} and \ref{resp-inj}.

\item  By \cite[Folgerung 2.2.24]{Abu2001} $_{A}C$ is locally Noetherian.
The result follows then from the isomorphism of categories $\mathrm{Rat}%
^{C}(_{A}\mathcal{M})\simeq \sigma \lbrack _{A}C]$ and \cite[27.3]{Wis88}.

\item  If $M\in \mathcal{M}^{C}$ is coflat, then $M_{R}$ is flat (by Remark
\ref{ko-flat}). Assume now that $_{A^{e}}A$ is projective. If $M_{R}$ is
flat, then by Proposition \ref{cot=Hom(ot)}
\begin{equation*}
M\square _{C}-\simeq \mathrm{Hom}_{A^{e}-}(A,-)\circ (M\otimes _{R}-)
\end{equation*}
is exact, i.e. $M$ is coflat.$\blacksquare $
\end{enumerate}
\end{Beweis}

\begin{corollary}
Let $(A,C)\in \mathcal{P}_{m}^{\alpha }.\;$If $R$ is semisimple \emph{(}e.g.
a field\emph{)}, then:

\begin{enumerate}
\item  $M\in \mathrm{Rat}^{C}(_{A}\mathcal{M})$ is injective if and only if $%
_{A}M$ is a direct summand of $_{A}C^{(\Lambda )}$ for some index set $%
\Lambda .$

\item  If $A$ is separable, then $C$ is right semisimple (i.e. every right $C
$-comodules is injective).
\end{enumerate}
\end{corollary}

\section{Coinduction Functors in $\mathcal{P}_{m}^{\protect\alpha }$}

By his study of the induced representations of quantum groups, Z. Lin (
\cite[3.2]{Lin93}, \cite{Lin94}) considered \emph{induction functors} for
\emph{admissible Hopf }$R$\emph{-pairings }over Dedekind rings. His aspect
was inspired by the induction functors in the theory of affine algebraic
groups and quantum groups. We generalize his results to the coinduction
functor for the category of \emph{measuring }$\alpha $\emph{-pairings} $%
\mathcal{P}_{m}^{\alpha }\subset \mathcal{P}_{m}$ and show that it is
isomorphic to a coinduction functor between categories of Type $\sigma
\lbrack M].$ Moreover we get as nice description of it as a composition of a
suitable $\mathrm{Hom}$-functor and\textrm{\ }a\textrm{\ }$\mathrm{Trace}$%
-functor.

\begin{punto}
\label{CoindMN}Let $A,B$ be $R$-algebras and $\xi :A\rightarrow B$ be an $R$%
-algebra morphism. Then every left $B$-module becomes a left $A$-module in a
canonical way and we get the so called \emph{restriction functor }$(-)_{\xi
}:$ $_{B}\mathcal{M}\rightarrow $ $_{A}\mathcal{M}.$ Considering $B$ with
the canonical $A$-bimodule structure, we have the functor $\mathrm{Hom}%
_{A-}(B,-):$ $_{A}\mathcal{M}\rightarrow $ $_{B}\mathcal{M}.$ Moreover $%
((-)_{\xi },\mathrm{Hom}_{A-}(B,-))$ is an adjoint pair of \emph{covariant}
functors through the functorial canonical isomorphisms
\begin{equation*}
\mathrm{Hom}_{A-}(M_{\xi },N)\simeq \mathrm{Hom}_{A-}(B\otimes
_{B}M,N)\simeq \mathrm{Hom}_{B-}(M,\mathrm{Hom}_{A-}(B,N)).
\end{equation*}
If we consider the \emph{induction functor }$B\otimes _{A}-:\mathcal{\ }_{A}%
\mathcal{M}\rightarrow $ $_{B}\mathcal{M},$ then $(B\otimes _{A}-,(-)_{\xi })
$ is an adjoint pair of \emph{covariant} functors through the functorial
canonical isomorphisms
\begin{equation*}
\mathrm{Hom}{\normalsize _{B-}}({\normalsize B}\otimes _{A}N,M)\simeq
\mathrm{Hom}_{A-}(N,\mathrm{Hom}{\normalsize _{B-}}(B,M))\simeq \mathrm{Hom}%
_{A-}(N,M_{\xi }).
\end{equation*}
\end{punto}

\begin{punto}
\textbf{The general coinduction functor}$.$\label{CoindKL} Let $A,B$ be $R$%
-Algebras and $\xi :A\rightarrow B$ be an $R$-algebra morphism. If $L$ is a
left $B$-module, then we get the covariant functor
\begin{equation}
\mathrm{HOM}_{A-}(B,-):={\normalsize \mathrm{Sp}(\sigma \lbrack }_{B}L%
{\normalsize ],}\mathrm{Hom}_{A-}{\normalsize (B,-)):}_{A}\mathcal{M}%
\rightarrow {\normalsize \sigma \lbrack }_{B}L{\normalsize ].}
\label{Big-Hom-sg}
\end{equation}
For every left $A$-module $K$ (\ref{Big-Hom-sg}) restricts to the covariant
\emph{coinduction functor}
\begin{equation}
\mathrm{Coind}_{K}^{L}{\normalsize (-):=\mathrm{Sp}(\sigma \lbrack }_{B}L%
{\normalsize ],}\mathrm{Hom}_{A-}{\normalsize (B,-)):\sigma \lbrack }_{A}K%
{\normalsize ]\rightarrow \sigma \lbrack }_{B}L{\normalsize ],}  \label{spMN}
\end{equation}
i.e. $\mathrm{Coind}_{K}^{L}(-)$ is defined through the commutativity of the
following diagram:
\begin{equation*}
\xymatrix{ {}_{A} {\mathcal M} \ar[rrrdd]^{{\rm HOM}_{A-} (B,-)} \ar[rrr]^{
{\rm Hom}_{A-} (B,-)} & & & {}_{B} {\mathcal M} \ar[dd]^{ {\rm Sp}
(\sigma[{}_{B}L],-) }\\ & & & \\ \sigma[{}_{A} K] \ar@{.>}[rrr]_{ {\rm
Coind} _K ^L (-)} \ar@{^{(}->}[uu] & & & \sigma[{}_{B} L] }
\end{equation*}
If $(L)_{\xi }$ is $K$-subgenerated as a left $A$-module, then $(-)_{\xi }:$
$_{B}\mathcal{M}\rightarrow $ $_{A}\mathcal{M}$ restricts to $(-)_{\xi }:%
\sigma \lbrack _{B}L]\rightarrow \sigma \lbrack _{A}K]$ and $((-)_{\xi },%
\mathrm{Coind}_{K}^{L}(-))$ turns to be an adjoint pair of covariant
functors.
\end{punto}

\begin{punto}
\label{th-box}\textbf{The ad-corestriction functor.}\emph{\ }Let $C,D$ be $R$%
-coalgebras and $\theta :D\rightarrow C$ be an $R$-coalgebra morphism. Then
we get the covariant \emph{corestriction functor}
\begin{equation}
(-)^{\theta }:\mathcal{M}^{D}\rightarrow \mathcal{M}^{C},\text{ }(M,\varrho
_{M})\mapsto (M,(id_{M}\otimes \theta )\circ \varrho _{M}).  \label{^th}
\end{equation}
On the other hand, consider $D$ as a left $C$-comodule through
\begin{equation*}
\varrho _{D}^{C}:D\overset{\Delta _{D}}{\longrightarrow }D\otimes _{R}D%
\overset{\theta \otimes id}{\longrightarrow }C\otimes _{R}D.
\end{equation*}
If $_{R}D$ is \emph{flat}, then for every $M\in \mathcal{M}^{C}$ the
cotensor product $M\square _{C}D$ becomes a right $D$-comodule through
\begin{equation*}
M\square _{C}D\overset{id\square _{C}\Delta _{D}}{\longrightarrow }M\square
_{C}(D\otimes _{R}D)\simeq (M\square _{C}D)\otimes _{R}D
\end{equation*}
and we get the \emph{ad-corestriction functor}
\begin{equation*}
-\square _{C}D:\mathcal{M}{\normalsize ^{C}\rightarrow \mathcal{M}^{D},}%
\text{ }M\mapsto M\square _{C}D.
\end{equation*}
\end{punto}

\begin{punto}
Let $Q=(B,D)\in \mathcal{P}_{m}^{\alpha }.$ For every $R$-algebra $A$ with $R
$-algebra morphism $\xi :A\rightarrow B$ we have the covariant functor
\begin{equation*}
\mathrm{HOM}_{A-}(B,-):=\mathrm{Rat}^{D}(-)\circ \mathrm{Hom}_{A-}(B,-):%
\text{ }_{A}\mathcal{M}\rightarrow {\normalsize \mathrm{Rat}^{D}(_{B}%
\mathcal{M})}.
\end{equation*}
If $P=(A,C)\in \mathcal{P}_{m}^{\alpha },$ then $\mathrm{HOM}_{A-}(B,-)$
restricts to the \emph{coinduction functor}\textbf{\ }from $P$ to $Q:$%
\begin{equation*}
\mathrm{Coind}_{P}^{Q}(-):\mathrm{Rat}^{C}(_{A}\mathcal{M}{\normalsize %
)\rightarrow \mathrm{Rat}^{D}(_{B}\mathcal{M}),}\text{ }M\mapsto \mathrm{Rat}%
^{D}(\mathrm{Hom}_{A-}(B,M)),
\end{equation*}
i.e. $\mathrm{Coind}_{P}^{Q}(-)$ is defined through the commutativity of the
following diagram:
\begin{equation*}
\xymatrix{ {}_{A} {\mathcal M} \ar[rrrdd]^{{\rm HOM}_{A-} (B,-)} \ar[rrr]^{
{\rm Hom}_{A-} (B,-)} & & & {}_{B} {\mathcal M} \ar[dd]^{ {\rm Rat}^{D} (-)}
\\ & & & \\ {\rm Rat}^{C}(_{A} {\mathcal M}) \ar@{.>}[rrr]_{ {\rm Coind} _P
^Q (-)} \ar@{^{(}->}[uu] & & & {\rm Rat}^{D}(_{B} {\mathcal M}) }
\end{equation*}
\end{punto}

\begin{proposition}
\label{adj-th-kp}Let $P=(A,C),$ $Q=(B,D)\in \mathcal{P}_{m}$ and $(\xi
,\theta ):(B,D)\rightarrow (A,C)$ be a morphism in $\mathcal{P}_{m}.$

\begin{enumerate}
\item  If $_{R}D$ is flat, then $((-)^{\theta },-\square _{C}D)$ is an
adjoint pair of covariant functors.

\item  If $P,Q\in \mathcal{P}_{m}^{\alpha }$ and $B$ is \emph{commutative},
then we have for every $N\in $ $_{A}\mathcal{M}:$%
\begin{equation*}
\mathrm{HOM}_{A-}(B,N)=\mathrm{HOM}_{A-}(B,\mathrm{Rat}^{C}(_{A}N)).
\end{equation*}
\end{enumerate}
\end{proposition}

\begin{Beweis}
\begin{enumerate}
\item  One can show easily that the mapping
\begin{equation*}
\Phi _{N,L}:\mathrm{Hom}^{{\normalsize D}}(N,L\square _{C}{\normalsize D}%
)\rightarrow \mathrm{Hom}^{C}(N^{\theta },L),\text{ }f\mapsto (id_{L}\square
_{C}\theta )\circ f
\end{equation*}
is an isomorphism with inverse $g\mapsto (g\square _{C}id_{D})\circ \varrho
_{N}$ and moreover that it is functorial in $N\in \mathcal{M}^{D}$ and $L\in
\mathcal{M}^{C}.$

\item  If $g\in \mathrm{HOM}_{A-}(B,N),$ then we have for all $a\in A$ and $%
b\in B:$%
\begin{equation*}
\begin{tabular}{lllll}
$a(g(b))$ & $=$ & $g(a\rightharpoondown b)$ & $=$ & $g(\xi (a)b)$ \\
& $=$ & $g(b\xi (a))$ & $=$ & $(\xi (a)g)(b)$ \\
& $=$ & $\sum g_{<0>}(b)<\xi (a),g_{<1>}>$ & $=$ & $\sum g_{<0>}(b)<a,\theta
(g_{<1>})>.$%
\end{tabular}
\end{equation*}
Consequently $g(B)\subseteq \mathrm{Rat}^{C}(_{A}N)$ and the result follows.$%
\blacksquare $
\end{enumerate}
\end{Beweis}

\begin{punto}
\label{funk-eq}Let $P=(A,C),$ $Q=(B,D)\in \mathcal{P}_{m}^{\alpha },$ $(\xi
,\theta ):(B,D)\rightarrow (A,C)$ be a morphism in $\mathcal{P}_{m}^{\alpha }
$ and denote the restriction of $(-)_{\xi }:$ $_{B}\mathcal{M}\rightarrow $ $%
_{A}\mathcal{M}$ on $\mathrm{Rat}^{D}(_{B}\mathcal{M})=\sigma \lbrack _{B}D]$
also with $(-)_{\xi }.$ Through the isomorphism of categories $\mathcal{M}%
^{C}\simeq \mathrm{Rat}^{C}(_{A}\mathcal{M})=\sigma \lbrack _{A}C]$ and $%
\mathcal{M}^{D}\simeq \mathrm{Rat}^{D}(_{B}\mathcal{M})=\sigma \lbrack _{B}D]
$ (compare Theorem \ref{cor-dicht}) we get an equivalence of functors $%
(-)^{\theta }\approx (-)_{\xi }.$ Considering the covariant functors (\ref
{cov-fund}) we get a commutative diagram of pairwise adjoint covariant
functors
\begin{equation}
\xymatrix{{\cal M}^{D} \ar[rrr]_{_{B} (-)} \ar@<-.3cm>[ddd]_{(-)^{\theta}} &
& & {\rm Rat} ^{D} (_{B} {\cal M}) \ar@<-.4cm>[lll]_{(-)^{D}}
\ar@<-.2cm>@{-}[rr] \ar@<.2cm>@{-}[rr] \ar@<-.3cm>[ddd]_{(-)_{\xi }} & &
\sigma [{}_{B} D] \ar@{^{(}->}[rrr]_{\iota _D} \ar@<-.3cm>[ddd]_{(-)_{\xi }}
& & & {}_{B} {\cal M} \ar@<-.4cm>[lll]_{{\rm{Sp}}(\sigma [_{B}D],-)}
\ar@<-.3cm>[ddd]_{(-)_{\xi }} \\ \\ \\ {\cal M}^C \ar[rrr]^{_{A} (-)}
\ar[uuu]_{- \Box _C D} & & & {\rm Rat} ^{C} (_{A} {\cal M})
\ar[uuu]_{\rm{Coind}_{P}^{Q}(-)} \ar@<-.2cm>@{-}[rr] \ar@<.2cm>@{-}[rr]
\ar@<.4cm>[lll]^{(-)^{C}} & & \sigma [{}_A C] \ar@{^{(}->}[rrr]^{\iota _C}
\ar[uuu]_{\rm{Coind}_{C} ^{D} (-)} & & & {}_{A} {\cal M}
\ar@<.4cm>[lll]^{{\rm{Sp}}(\sigma [ _{A}C],-)} \ar[uuu]_{{\rm Hom}_{A-}
(B,-)} }  \label{adj-diag}
\end{equation}
\end{punto}

\begin{theorem}
\label{phi-ad}Let $P=(A,C),$ $Q=(B,D)\in \mathcal{P}_{m}^{\alpha }$ \emph{(}%
so that in particular $_{R}C$ and $_{R}D$ are flat\emph{)} and $(\xi ,\theta
):(B,D)\rightarrow (A,C)$ be a morphism in $\mathcal{P}_{m}^{\alpha }.$
Through the isomorphisms of categories $\mathcal{M}^{C}\simeq \mathrm{Rat}%
^{C}(_{A}\mathcal{M})=\sigma \lbrack _{A}C]$ and $\mathcal{M}^{D}\simeq
\mathrm{Rat}^{D}(_{B}\mathcal{M})=\sigma \lbrack _{B}D]$ \emph{(}compare
Theorem \ref{equal}\emph{)} the following functors are equivalent
\begin{equation*}
\begin{tabular}{lllll}
$-\square _{C}D$ & $:$ & $\mathcal{M}^{C}$ & $\rightarrow $ & $\mathcal{M}%
^{D},$ \\
$\mathrm{Coind}_{P}^{Q}(-)$ & $:$ & $\mathrm{Rat}^{C}(_{A}\mathcal{M})$ & $%
\rightarrow $ & $\mathrm{Rat}^{D}(_{B}\mathcal{M}),$ \\
$\mathrm{Hom}_{A^{e}-}(A,-\otimes _{R}D)$ & $:$ & $\mathrm{Rat}^{C}(_{A}%
\mathcal{M})$ & $\rightarrow $ & $\mathrm{Rat}^{D}(_{B}\mathcal{M}),$ \\
$\mathrm{Coind}_{C}^{D}(-)$ & $:$ & $\sigma \lbrack _{A}C]$ & $\rightarrow $
& $\sigma \lbrack _{B}D].$%
\end{tabular}
\end{equation*}
\end{theorem}

\begin{Beweis}
Consider for every $N\in \mathcal{M}^{C}$ the \emph{injective }$R$-linear
mapping
\begin{equation*}
\gamma _{N}:=(\alpha _{N}^{Q})|_{N\square _{C}{\normalsize D}}:N\square _{C}%
{\normalsize D}\rightarrow \mathrm{Hom}_{R}({\normalsize B},N),\text{ }\sum
n_{i}\otimes d_{i}\mapsto \lbrack b\mapsto \sum n_{i}<b,d_{i}>].
\end{equation*}
Then we have for all $a\in A$ and $b\in B:$%
\begin{equation*}
\begin{tabular}{lllll}
$\gamma _{N}(\sum n_{i}\otimes d_{i})(a\rightharpoondown b)$ & $=$ & $\sum
n_{i}<a\rightharpoondown b,d_{i}>$ &  &  \\
& $=$ & $\sum n_{i}<b,d_{i}\leftharpoonup a>$ &  &  \\
& $=$ & $\gamma _{N}(\sum n_{i}\otimes d_{i}\leftharpoonup a)(b)$ &  &  \\
& $=$ & $\gamma _{N}(\sum an_{i}\otimes d_{i})(b)$ & (compare Lemma \ref
{cot=Hom(ot)} (1)) &  \\
& $=$ & $\sum an_{i}<b,d_{i}>$ &  &  \\
& $=$ & $a(\gamma _{N}(n_{i}\otimes d_{i})(b)),$ &  &
\end{tabular}
\end{equation*}
i.e. $\gamma _{N}(N\square _{C}D)\subset \mathrm{Hom}_{A-}(B,N).$ Moreover
we have for arbitrary $\sum n_{i}\otimes d_{i}\in N\square _{C}D$ and $b,%
\widetilde{b}\in B:$%
\begin{equation*}
\begin{tabular}{lllll}
$\gamma _{N}(\widetilde{b}(\sum n_{i}\otimes d_{i}))(b)$ & $=$ & $\gamma
_{N}(\sum n_{i}\otimes \widetilde{b}\rightharpoonup d_{i})(b)$ & $=$ & $\sum
n_{i}<b,\widetilde{b}\rightharpoonup d_{i}>$ \\
& $=$ & $\sum n_{i}<b\widetilde{b},d_{i}>$ & $=$ & $\gamma _{N}(\sum
n_{i}\otimes d_{i})(b\widetilde{b})$ \\
& $=$ & $(\widetilde{b}\gamma _{N}(\sum n_{i}\otimes d_{i}))(b),$ &  &
\end{tabular}
\end{equation*}
i.e. $\gamma _{N}$ is $B$-linear. But $_{R}D$ is flat, so $N\square _{C}D\in
\mathcal{M}^{D}$ by Corollary \ref{cot-com} and it follows by Lemma \ref
{clos} that $\gamma _{N}(N\square _{C}D)\subset \mathrm{HOM}_{A-}(B,N).$ Now
we show that the following $R$-linear mapping is well defined
\begin{equation*}
\beta _{N}:\mathrm{HOM}{\normalsize _{A-}(B,N)\rightarrow N\square _{C}D,}%
\text{ }f\mapsto \sum f_{<0>}(1_{B})\otimes f_{<1>}.
\end{equation*}
For all $f\in \mathrm{HOM}_{A-}(B,N),$ $a\in A$ and $b\in B$ we have
\begin{equation*}
\begin{tabular}{lll}
$\gamma _{N}(\sum a(f_{<0>}(1_{B}))\otimes f_{<1>})(b)$ & $=$ & $\sum
a(f_{<0>}(1_{B}))<b,f_{<1>}>$ \\
& $=$ & $\sum f_{<0>}(a\rightharpoondown 1_{B}))<b,f_{<1>}>$ \\
& $=$ & $\sum f_{<0>}(\xi (a))<b,f_{<1>}>$ \\
& $=$ & $\sum (\xi (a)f_{<0>})(1_{B})<b,f_{<1>}>$ \\
& $=$ & $\sum f_{<0><0>}(1_{B})<\xi (a),f_{<0><1>}><b,f_{<1>}>$ \\
& $=$ & $\sum f_{<0>}(1_{B})<\xi (a),f_{<1>1}><b,f_{<1>2}>$ \\
& $=$ & $\sum f_{<0>}(1_{B})<\xi (a)b,f_{<1>}>$ \\
& $=$ & $\sum f_{<0>}(1_{B})<a\rightharpoondown b,f_{<1>}>$ \\
& $=$ & $\sum f_{<0>}(1_{B})<b,f_{<1>}\leftharpoonup a>$ \\
& $=$ & $\gamma _{N}(\sum f_{<0>}(1_{B})\otimes f_{<1>}\leftharpoonup a)(b),$%
\end{tabular}
\end{equation*}
i.e. $\sum a(f_{<0>}(1_{B}))\otimes f_{<1>}=\sum f_{<0>}(1_{B})\otimes
f_{<1>}\leftharpoonup a$ (since $\gamma _{N}$ is injective). It follows then
by Proposition \ref{cot=Hom(ot)} (1) that $\sum f_{<0>}(1_{B})\otimes
f_{<1>}\in N\square _{C}D,$ i.e. $\beta _{N}$ is well defined. Moreover, we
have for all $f\in \mathrm{HOM}_{A-}(B,N)$ and $b\in B:$%
\begin{equation*}
\begin{tabular}{lll}
$(\gamma _{N}\circ \beta _{N})(f)(b)$ & $=$ & $\gamma _{N}(\sum
f_{<0>}(1_{B})\otimes f_{<1>})(b)$ \\
& $=$ & $\sum f_{<0>}(1_{B})<b,f_{<1>}>$ \\
& $=$ & $(bf)(1_{B})=f(b),$%
\end{tabular}
\end{equation*}
hence $\gamma _{N}\circ \beta _{N}=id.$ Obviously $\beta _{N}\circ \gamma
_{N}=id.$ Consequently $\gamma _{N}$ and $\beta _{N}$ are isomorphisms. It
is easy to show that $\gamma _{N}$ and $\beta _{N}$ are functorial in $N,$
hence $-\square _{C}D\approx \mathrm{Coind}_{P}^{Q}(-).$ The equivalences $%
\mathrm{Coind}_{P}^{Q}(-)\approx \mathrm{Coind}_{C}^{D}(-)$ and $-\square
_{C}D\approx \mathrm{Hom}_{A^{e}-}(A,-\otimes _{R}D)$ follow now by Theorem
\ref{cor-dicht} and Proposition \ref{cot=Hom(ot)} (2), respectively.$%
\blacksquare $
\end{Beweis}

\begin{punto}
\label{otC=rtc}Let $Q=(B,D)\in \mathcal{P}_{m}^{\alpha }$ and consider the
trivial $R$-pairing $P=(R,R)\in \mathcal{P}_{m}^{\alpha }$ with the morphism
of measuring $R$-pairings$\ (\eta _{B},\varepsilon _{D}):(B,D)\rightarrow
(R,R).$ Then we have for every $M\in \mathcal{M}^{R}\simeq \mathcal{M}_{R}$%
\begin{equation*}
\mathrm{Coind}_{P}^{Q}(-):=\mathrm{HOM}_{R}(B,-)\simeq -\otimes _{R}D.
\end{equation*}
Notice that $\digamma \simeq (-)^{\varepsilon }:\mathcal{M}^{D}\rightarrow
\mathcal{M}_{R},$ where $\digamma $ is the forgetful functor, hence $%
(\digamma ,\mathrm{Coind}_{P}^{Q}(-))$ is an adjoint pair of covariant
functors.
\end{punto}

\begin{punto}
\textbf{Universal Property. }Let $P=(A,C),$ $Q=(B,D)\in \mathcal{P}%
_{m}^{\alpha }$ and $(\xi ,\theta ):(B,D)\rightarrow (A,C)$ be a morphism in
$\mathcal{P}_{m}^{\alpha }.$ Then $\mathrm{Coind}_{P}^{Q}(-)$ has the
following universal property: if $N\in \mathcal{M}^{D},$ $M\in \mathcal{M}%
^{C}$ and $\phi \in \mathrm{Hom}^{C}(N^{\theta },M),$ then there exists a
unique $\widetilde{\phi }\in \mathrm{Hom}^{D}(N,\mathrm{Coind}_{P}^{Q}(M)),$
such that $\phi (n)=\widetilde{\phi }(n)(1_{B})$ for every $n\in N.$
\end{punto}

In what follows we list some properties of the coinduction functor:

\begin{punto}
\label{ind-prop}Let $P=(A,C),$ $Q=(B,D)\in \mathcal{P}_{m}^{\alpha }$ and $%
(\xi ,\theta ):(B,D)\rightarrow (A,C)$ be a morphism in $\mathcal{P}%
_{m}^{\alpha }$.

\begin{enumerate}
\item  $\mathrm{Coind}_{P}^{Q}(-)$ respects direct limits: if $\{N_{\lambda
}\}_{\Lambda }$ is a directed system in $\mathrm{Rat}^{C}(_{A}\mathcal{M}),$
then
\begin{equation*}
\mathrm{Coind}{\normalsize _{P}^{Q}}(\underrightarrow{lim}N_{\lambda
})\simeq \underrightarrow{lim}N_{\lambda }\square _{C}D\simeq
\underrightarrow{lim}(N_{\lambda }\square _{C}D)=\underrightarrow{lim}%
\mathrm{Coind}_{P}^{Q}(N_{\lambda }).
\end{equation*}

\item  $\mathrm{Rat}^{D}(-)$ $\&$ $\mathrm{Hom}_{A-}(B,-)$ are left-exact,
hence
\begin{equation*}
\mathrm{Coind}_{P}^{Q}(-):=\mathrm{Rat}^{D}(-)\circ \mathrm{Hom}_{A-}(B,-)
\end{equation*}
is left-exact. If moreover $_{A}B$ is projective (hence $\mathrm{Hom}%
_{A-}(B,-)$ is exact) and $\mathrm{Rat}^{D}(-)$ is exact, then $\mathrm{Coind%
}_{P}^{Q}(-)$ is exact.

\item  $\mathrm{Coind}_{P}^{Q}(-)\simeq -\square _{C}D$ is exact if and only
if $D$ is coflat in $^{C}\mathcal{M}.$ If $R$ is a QF ring, then $\mathrm{%
Coind}_{P}^{Q}(-)$ is exact if and only if $D$ is injective in $^{C}\mathrm{%
Rat}(\mathcal{M}_{A}).$

\item  By Lemma \ref{adj-th-kp} (1) $((-)^{\theta },-\square _{C}D)$ is an
adjoint pair of covariant functors, hence $\mathrm{Coind}_{P}^{Q}(-)\simeq
-\square _{C}D$ respects inverse Limit, i.e. direct products, kernels and
injective objects (since $(-)^{\theta }:\mathcal{M}^{D}\rightarrow \mathcal{M%
}^{C}$ is exact). In particular, if $C$ is injective in $\mathrm{Rat}%
^{C}(_{A}\mathcal{M}),$ then $D\simeq C\square _{C}D\simeq \mathrm{Coind}%
_{P}^{Q}(C)$ is injective in $\mathrm{Rat}^{D}(_{B}\mathcal{M}).$

\item  Let $A$ be \emph{separable}. Then $-\square _{C}D\simeq \mathrm{Coind}%
_{P}^{Q}(-)\simeq $ $\mathrm{Hom}_{A^{e}-}(A,-\otimes _{R}D)$ is exact, i.e.
$D$ is coflat in $^{C}\mathcal{M}$. If moreover $R$ is a QF ring, then $D$
is injective in $^{C}\mathcal{M}.$
\end{enumerate}
\end{punto}

A version of the following result was obtained by Y. Doi \cite[Proposition 5]
{Doi81} in the case of a base field:

\begin{proposition}
Let $P=(A,C),$ $Q=(B,D)\in \mathcal{P}_{m}^{\alpha }$ and $(\xi ,\theta
):(B,D)\rightarrow (A,C)$ be a morphism in $\mathcal{P}_{m}^{\alpha }.$ If $%
R\;$is a QF ring, then the following are equivalent:

\begin{enumerate}
\item  The functor $\mathrm{Coind}_{P}^{Q}(-):$ $\mathrm{Rat}^{C}(_{A}%
\mathcal{M})\rightarrow \mathrm{Rat}^{D}(_{B}\mathcal{M})$ is exact;

\item  $D$ is coflat in $^{C}\mathcal{M};$

\item  $D$ is injective in $^{C}\mathrm{Rat}(\mathcal{M}_{A});$

\item  If $M$ is an injective left $D$-comodule that is flat in $\mathcal{M}%
_{R},$ then $M$ is injective in $^{C}\mathrm{Rat}(\mathcal{M}_{A}).$
\end{enumerate}
\end{proposition}

\begin{Beweis}
(1) $\Leftrightarrow $ (2) Follows from the isomorphism of functors $\mathrm{%
Coind}_{P}^{Q}(-)\simeq -\square _{C}D:\mathcal{M}^{C}\rightarrow \mathcal{M}%
^{D}.$

(2) $\Leftrightarrow $ (3) By Remark\ \ref{flat} $_{R}D$ is flat, so the
equivalence follows from Lemma \ref{inj=cof}.

(2) $\Rightarrow $ (4) Let $M$ be a left $D$-comodule and assume that $M$ is
injective in $^{D}\mathrm{Rat}(\mathcal{M}_{B})$ and flat in $\mathcal{M}%
_{R}.$ Then $M$ is coflat in $^{D}\mathcal{M}$ (by Lemma \ref{inj=cof}) and
we have (by Lemma \ref{koass-cot}) an isomorphism of functors
\begin{equation*}
-\square _{C}M\simeq (-\square _{C}D)\square _{D}M:\mathcal{M}%
^{C}\rightarrow \mathcal{M}_{R}
\end{equation*}
By assumption $-\square _{C}D:\mathcal{M}^{C}\rightarrow \mathcal{M}^{D}$
and $-\square _{D}M:\mathcal{M}^{D}\rightarrow \mathcal{M}_{R}$ are exact
and so $-\square _{C}M$ is exact. By Lemma \ref{inj=cof} $M$ is injective in
$^{C}\mathrm{Rat}(\mathcal{M}_{A}).$

(4) $\Rightarrow $ (3) Since $R$ is injective, $D$ is injective in $^{D}%
\mathrm{Rat}(\mathcal{M}_{B})\simeq $ $^{D}\mathcal{M}.$ It follows then
from the assumption that $D$ is injective in $^{C}\mathrm{Rat}(\mathcal{M}%
_{A}).\blacksquare $
\end{Beweis}

\qquad As a consequence of Theorem \ref{cor-dicht} and \cite[Proposition
3.23]{Abu2} we get:

\begin{corollary}
\label{ind-loc}Let $R$ be Noetherian, $A,B$ be $R$-algebras and $\xi
:A\rightarrow B$ be an $R$-algebra morphism.

\begin{enumerate}
\item  Let $A,B$ be cofinitely $R$-cogenerated $\alpha $-algebras, $%
P:=(A,A^{\circ }),$ $Q:=(B,B^{\circ })$ and consider the morphism of
measuring $\alpha $-pairings $(\xi ,\xi ^{\circ }):(B,B^{\circ })\rightarrow
(A,A^{\circ }).$ Then we have for every right $A^{\circ }$-comodule $N:$%
\begin{equation*}
\mathrm{Coind}_{P}^{Q}(N)=\{f\in \mathrm{Hom}_{A-}(B,N)|\text{ }Bf\text{ is
finitely generated in }\mathcal{M}_{R}\}.
\end{equation*}

\item  Let ${\frak{F}}_{A},$ ${\frak{F}}_{B}$ be cofinitely $R$-cogenerated $%
\alpha $-filters of $R$-cofinite $A$-ideals, $B$-ideals respectively and
consider $A,$ $B$ as a left linear topological $R$-algebra with the induced
left linear topologies $\frak{T}({\frak{F}}_{A}),$ $\frak{T}({\frak{F}}_{B})$
respectively. If $\xi :A\rightarrow B$ be an $R$-algebra morphism that is
continuous with respect to $(\frak{T}({\frak{F}}_{A}),\frak{T}({\frak{F}}%
_{B})),$ $P:=(A,A_{\frak{F}_{A}}^{\circ })$ and $Q:=(B,B_{\frak{F}%
_{B}}^{\circ }),$ then we have for every $N\in \mathcal{M}^{A_{\frak{F}%
_{A}}^{\circ }}:$%
\begin{equation*}
\mathrm{Coind}_{P}^{Q}(N)=\{f\in \mathrm{Hom}_{A-}(B,N)|\text{ }(0:f)\supset
\widetilde{I}\text{ for some }\widetilde{I}\in {\frak{F}}{\normalsize _{B}}%
\}.
\end{equation*}
\end{enumerate}
\end{corollary}

\section{Hopf $R$-pairings}

\begin{definition}
Let $H$ be an $R$-bialgebra. An $H$-ideal, which is also an $H$-coideal, is
called a \emph{bi-ideal.} If $H$ is a Hopf $R$-algebra with antipode $S_{H}$
and $J\subset H$ is an $H$-bi-ideal with $S_{H}(J)\subset J,$ then $J$ is
called a \emph{Hopf ideal}.
\end{definition}

\begin{punto}
\textbf{The category }$\mathcal{P}_{Big}.$ A \emph{bialgebra }$R$\emph{%
-pairing }is an $R$-pairing $P=(H,K),$ where $H,K$ are $R$-bialgebras and $%
\kappa _{P}:H\rightarrow K^{\ast },$ $\chi _{P}:K\rightarrow H^{\ast }$ are $%
R$-algebra morphisms. For bialgebra $R$-pairings $(H,K),(Y,K)$ a morphism of
$R$-pairings $(\xi ,\theta ):(Y,Z)\rightarrow (H,K)$ is said to be a \emph{%
morphism of bialgebra }$R$\emph{-pairings,} if $\xi :H\rightarrow Y$ and $%
\theta :Z\rightarrow K$ are $R$-bialgebra morphisms. With $\mathcal{P}%
_{Big}\subset \mathcal{P}_{m}$ we denote the subcategory of bialgebra $R$%
-pairings and with $\mathcal{P}_{Big}^{\alpha }\subset \mathcal{P}_{Big}$
the \emph{full} subcategory, whose objects satisfy the{\normalsize \ }$%
\alpha $-condition.

If $P=(H,K)\in \mathcal{P}_{Big},$ $Z\subset K$ is a (pure) $R$-subbialgebra
and $J\subset H$ is an $H$\emph{-bi-ideal} with $<J,Z>=0,$ then $Q=(H/J,Z)$
is a bialgebra $R$-pairing, $(\pi _{J},\iota _{Z}):(H/J,Z)\rightarrow (H,K)$
is a morphism in $\mathcal{P}_{Big}$ and $Q\subset P$ is called a (pure)
\emph{bialgebra }$R$\emph{-subpairing.} Obviously $\mathcal{P}_{Big}^{\alpha
}\subset \mathcal{P}_{Big}$ is closed under pure bialgebra $R$-subpairings.
\end{punto}

\begin{punto}
\textbf{The category }$\mathcal{P}_{Hopf}.$ A \emph{Hopf }$R$\emph{-pairing}
$P=(H,K)$ is a bialgebra $R$-pairing with $H,K$ Hopf $R$-algebras. With $%
\mathcal{P}_{Hopf}\subset \mathcal{P}_{Big}$ we denote the \emph{full }%
subcategory of Hopf $R$-pairings and with $\mathcal{P}_{Hopf}^{\alpha
}\subset \mathcal{P}_{Hopf}$ the \emph{full }subcategory, whose objects
satisfy the{\normalsize \ }$\alpha $-condition. If $P=(H,K)$ is a Hopf $R$%
-pairing, $Z\subset K$ a (pure) Hopf $R$-subalgebra and $J\subset H$ a \emph{%
Hopf ideal} with $<J,Z>=0,$ then $Q:=(H/J,Z)$ is a Hopf $R$-pairing, $(\pi
_{J},\iota _{Z}):(H/J,Z)\rightarrow (H,K)$ is a morphism in $\mathcal{P}%
_{Hopf}$ and $Q\subset P$ is called a (pure) \emph{Hopf }$R$\emph{-subpairing%
}. Obviously $\mathcal{P}_{Hopf}^{\alpha }\subset \mathcal{P}_{Hopf}$ is
closed under pure Hopf $R$-subpairings.
\end{punto}

\begin{example}
Let $R$ be Noetherian and $H$ be an $\alpha $-bialgebra (respectively a Hopf
$\alpha $-algebra). Then $H^{\circ }$ is by (\cite[Theorem 2.8]{AG-TW2000})
an $R$-bialgebra (respectively a Hopf $R$-algebra). Moreover it is easy to
see that $(H,H^{\circ })$ is a bialgebra $\alpha $-pairing (respectively a
Hopf $\alpha $-pairing).
\end{example}

\begin{remarks}
\begin{enumerate}
\item  (Compare \cite{Tak92}) If $P=(H,K)$ is a Hopf $R$-pairing, then
\begin{equation*}
<S_{H}(h),k>=<h,S_{K}(k)>\text{ for all }h\in H\text{ and }k\in K.
\end{equation*}

\item  Let $R$ be Noetherian. If $P=(H,K)$ is a bialgebra $R$-pairing
(respectively a Hopf $R$-pairing), then $\kappa _{P}(H)\subset K^{\circ }$
and $\chi _{P}(K)\subset H^{\circ }.$ If $(H,K)\in \mathcal{P}_{Big}$ and $%
H\in \mathbf{Big}_{R}^{\alpha }$ (respectively $K\in \mathbf{Big}%
_{R}^{\alpha }$), then $\chi _{P}:K\rightarrow H^{\circ }$ (respectively $%
\kappa _{P}:H\rightarrow K^{\circ }$) is an $R$-bialgebra morphism.
\end{enumerate}
\end{remarks}

\subsection*{Quasi-Admissible filters.}

\qquad

By the study of induced representations of quantum groups, Z. Lin \cite
{Lin93} and M. Takeuchi \cite{Tak94} studied the so called \emph{admissible
filters}{\normalsize \ }of ideals of a Hopf $R$-algebra over arbitrary
(Dedekind) rings. In what follows we introduce what we call the \emph{%
quasi-admissible filters} and generalize some of their results to the class
of (not necessarily cofinitary) \emph{quasi-admissible }$\alpha $\emph{%
-filters.}

\begin{punto}
\label{filt-2}Let $A,B$ be $R$-algebras and $\frak{F}_{A},$ $\frak{F}_{B}$
be filters consisting of $R$-cofinite $A$-ideals, $B$-ideals respectively.
Then the filter basis
\begin{equation}
\frak{F}_{A}\times \frak{F}_{B}:=\{\mathrm{Im}(\iota _{I}\otimes id_{B})+%
\mathrm{Im}(id_{A}\otimes \iota _{J})|\text{ }I\in \frak{F}_{A},\text{ }J\in
\frak{F}_{B}\}  \label{filt-mal}
\end{equation}
induces on $A\otimes _{R}B$ a topology $\frak{T}(\frak{F}_{A}\times \frak{F}%
_{B}),$ such that $(A\otimes _{R}B,\frak{T}(\frak{F}_{A}\times \frak{F}_{B}))
$ is a linear topological $R$-algebra and $\frak{F}_{A}\times \frak{F}_{B}\ $%
is a neighbourhood basis of $0_{A\otimes _{R}B}.$
\end{punto}

\begin{punto}
\label{zul-filt}{\normalsize \ }Let $H$ be an $R$-bialgebra (that is not a
Hopf $R$-algebra), $\frak{F}\subset \mathcal{K}_{H}$ be a filter and
consider the induced linear topological $R$-algebras $(H,\frak{T}(\frak{F}))$
and $(H\otimes _{R}H,\frak{T}(\frak{F}\times \frak{F})).$ We call $\frak{F}$
\emph{quasi-admissible,} if $\Delta _{H}:H\rightarrow H\otimes _{R}H$ and $%
\varepsilon _{H}:H\rightarrow R$ are continuous, i.e. if $\frak{F}$
satisfies the following axioms:
\begin{equation}
\begin{tabular}{ll}
(A1) & $\forall $ $I,J\in \frak{F}$ there exists $L\in \frak{F},$ such that $%
\Delta _{H}(L)\subseteq \mathrm{Im}(\iota _{I}\otimes id_{H})+\mathrm{Im}%
(id_{H}\otimes \iota _{J})$%
\end{tabular}
\label{A1}
\end{equation}
and
\begin{equation}
\begin{tabular}{ll}
(A2) & $\exists $ $I\in \frak{F},$ such that $\mathrm{Ker}(\varepsilon
_{H})\supset I.$%
\end{tabular}
\label{A2}
\end{equation}
If $H$ is a Hopf $R$-algebra, then we call a filter{\normalsize \ }$\frak{F}%
\subset \mathcal{K}_{H}\ $\emph{quasi-admissible}, if it satisfies (A1),
(A2) as well as
\begin{equation}
\begin{tabular}{ll}
(A3) & $\text{for every }I\in \frak{F}\text{ there exists }J\in \frak{F},%
\text{ such that }S_{H}(J)\subseteq I$%
\end{tabular}
\label{A3}
\end{equation}
(i.e. if $\Delta _{H},\varepsilon _{H}$ and $S_{H}$ are continuous). In \cite
{Lin93} and \cite{Tak94}, a cofinitary quasi-admissible filter of $R$%
-cofinite $H$-ideals (for a Hopf $R$-algebra $H$) is called \emph{admissible}%
.
\end{punto}

\begin{definition}
We call an $R$-bialgebra (respectively Hopf $R$-algebra) $H$ a \emph{%
quasi-admissible }$R$\emph{-bialgebra} (respectively a \emph{%
quasi-admissible Hopf }$R$\emph{-algebra}), if the class of $R$-cofinite $H$%
-ideals $\mathcal{K}_{H}$ is a \emph{quasi-admissible} filter.
\end{definition}

\begin{lemma}
\label{noeth-zul}If the ground ring $R$ is Noetherian, then every $R$%
-bialgebra \emph{(}Hopf $R$-algebra\emph{)} is quasi-admissible.
\end{lemma}

\begin{Beweis}
Let $H$ be an $R$-bialgebra. Since $R$ is Noetherian, $\mathcal{K}_{H}$ is a
filter. Moreover, $H\simeq R\oplus \mathrm{Ker}(\varepsilon _{H}),$ hence $%
\mathrm{Ker}(\varepsilon _{H})\in \mathcal{K}_{H}.$ Let $I,J\in \mathcal{K}%
_{H}$ and set $L:=\mathrm{Im}(\iota _{I}\otimes id_{H})+\mathrm{Im}%
(id_{H}\otimes \iota _{J}).$ Notice that $(H\otimes _{R}H)/L\simeq
H/I\otimes _{R}H/J$ (e.g. \cite[II-3.6, III-4.2]{Bou74}), hence $L\in
\mathcal{K}_{H\otimes _{R}H}.$ By definition $\Delta :H\rightarrow H\otimes
_{R}H$ is an $R$-algebra morphism and it follows, by the assumption $R$ is
Noetherian, that $\Delta ^{-1}(L)\vartriangleleft H$ is an $R$-cofinite
ideal. Consequently $H$ is a quasi-admissible $R$-bialgebra.

If $H$ is moreover a Hopf $R$-algebra, then $S_{H}:H\rightarrow H$ is an $R$%
-algebra antimorphism and it follows, from the assumption $R$ is Noetherian,
that for every $R$-cofinite ideal $I\vartriangleleft H$ the $H$-ideal $%
S_{H}^{-1}(I)\vartriangleleft H$ is $R$-cofinite. Consequently $H$ is a
quasi-admissible Hopf $R$-algebra.$\blacksquare $
\end{Beweis}

\begin{definition}
(\cite{Tak81}) An $R$-coalgebra $C$ is called \emph{infinitesimal flat,}%
{\normalsize \ }if $C=\underrightarrow{lim}C_{\lambda }$ for a directed
system of \emph{finitely generated projective} $R$-subcoalgebras $%
\{C_{\lambda }\}_{\Lambda }.$
\end{definition}

\begin{proposition}
\label{K(H)-zul}Let $H$ be an $R$-bialgebra \emph{(}respectively a Hopf $R$%
-algebra\emph{)} and $\frak{F}\subset \mathcal{K}_{H}$ be\ a
quasi-admissible filter.

\begin{enumerate}
\item  If $R$ is Noetherian and $\frak{F}$ is an $\alpha $-filter, then $H_{%
\frak{F}}^{\circ }$ is an $R$-bialgebra \emph{(}respectively a Hopf $R$%
-algebra\emph{)} and $(H,H_{\frak{F}}^{\circ })$ is a bialgebra $\alpha $%
-pairing \emph{(}respectively a Hopf $\alpha $-pairing\emph{)}.

\item  If $\frak{F}$\ is moreover cofinitary, then $H_{\frak{F}}^{\circ }$
is an \emph{infinitesimal flat} $R$-bialgebra \emph{(}Hopf $R$-algebra\emph{)%
} and $(H,H_{\frak{F}}^{\circ })$ is a bialgebra $\alpha $-pairing \emph{(}a
Hopf $\alpha $-pairing\emph{)}.
\end{enumerate}
\end{proposition}

\begin{Beweis}
\begin{enumerate}
\item  Let $H$ be an $R$-bialgebra. Obviously $H_{\frak{F}}^{\circ }\subset
H^{\circ }$ is an $H$-subbimodule under the regular left and the right $H$%
-actions (\ref{regular}) and so an $R$-coalgebra by Theorem \ref{R(G)-co}.
If $f(I)=0$ and $g(J)=0$ for $I,J\in \frak{F},$ then there exists by (\ref
{A1}) some $L\in \frak{F},$ such that $\Delta (L)\subseteq \mathrm{Im}(\iota
_{I}\otimes id_{H})+\mathrm{Im}(id_{H}\otimes \iota _{J}).$ Consequently $%
\Delta ^{\circ }(f\otimes g)(L)=(f\underline{\otimes }g)(\Delta (L))=0,$
i.e. $f\star g\in \mathrm{An}(L)\subset H_{\frak{F}}^{\circ }.$ By (\ref{A2}%
) $\varepsilon _{H}\in H_{\frak{F}}^{\circ }$ and so $H_{\frak{F}}^{\circ
}\subset H^{\ast }$ is an $R$-subalgebra. It is easy to see that $\Delta
^{\circ }:H_{\frak{F}}^{\circ }\otimes _{R}H_{\frak{F}}^{\circ }\rightarrow
H_{\frak{F}}^{\circ }$ and $\varepsilon ^{\circ }:R\rightarrow H_{\frak{F}%
}^{\circ }$ are coalgebra morphisms, i.e. $H_{\frak{F}}^{\circ }$ is an $R$%
-bialgebra. If $H\;$is a Hopf $R$-algebra with Antipode $S,$ then it follows
from (\ref{A3}) that $S^{\circ }(H_{\frak{F}}^{\circ })\subseteq H_{\frak{F}%
}^{\circ },$ hence $H_{\frak{F}}^{\circ }$ is a Hopf $R$-algebra with
antipode $S^{\circ }.$

\item  See \cite{Tak94}.$\blacksquare $
\end{enumerate}
\end{Beweis}

\qquad As a consequence of Lemma \ref{noeth-zul} and Proposition \ref
{K(H)-zul} we get

\begin{corollary}
\label{zul-noeth}Let $R$ be Noetherian. If $H$ is an $\alpha $-bialgebra
\emph{(}respectively a Hopf $\alpha $-algebra\emph{)}, then $H^{\circ }$ is
an $R$-bialgebra \emph{(}respectively a Hopf $R$-algebra\emph{)}. If $H$ is
cofinitary, then $H^{\circ }$ is an infinitesimal flat $R$-bialgebra \emph{(}%
respectively Hopf $R$-algebra\emph{)}.
\end{corollary}

\begin{proposition}
\label{dicht-zul}Let $H$ be an $R$-bialgebra, $\frak{F}${\normalsize \ }be a
quasi-admissible filter of $R$-cofinite $H$-ideals and consider $H$ as a
left linear topological $R$-algebra with the induced left linear topology $%
\frak{T}(\frak{F}).$ If $R$ is an injective cogenerator, then the following
are equivalent:

\begin{enumerate}
\item  $\frak{T}(\frak{F})$ is Hausdorff;

\item  the canonical $R$-linear mapping $\lambda :H\rightarrow H_{\frak{F}%
}^{\circ \ast }$ is injective;

\item  $H_{\frak{F}}^{\circ }\subset H^{\ast }$ is dense;

\item  $\sigma \lbrack _{H_{\frak{F}}^{\circ }}H]=\sigma \lbrack _{H^{\ast
}}H].$
\end{enumerate}
\end{proposition}

\begin{Beweis}
By assumption $H/I$ is $R$-cogenerated for every $I\in \frak{F}$ (hence $I=%
\mathrm{KeAn}(I)$ by \cite[28.1]{Wis88}) and so
\begin{equation*}
\overline{0_{A}}:=\bigcap_{I\in \frak{F}}I=\bigcap_{I\in \frak{F}}\mathrm{%
KeAn}(I)=\mathrm{Ke}(\sum\limits_{I\in \frak{F}}\mathrm{An}(I))=\mathrm{Ke}%
(H_{\frak{F}}^{\circ })=\mathrm{Ker}(\lambda ).
\end{equation*}
Since $R$ is an injective cogenerator, the equivalence (2) $\Leftrightarrow $
(3)\ follows from \cite[Theorem 1.8 (2)]{Abu1}. By assumption $\frak{F}\ $is
quasi-admissible, hence $H_{\frak{F}}^{\circ }\subset H^{\ast }$ is an $R$%
-subalgebra and the equivalence (3)\ $\Leftrightarrow $ (4) follows by Lemma
\ref{dense}.$\blacksquare $
\end{Beweis}

\qquad The proof of the following Proposition is along the lines of the
proof of \cite[Theorem 4.10]{AG-TL2001}:

\begin{proposition}
\label{pro-2H}Let $H,K$ be $R$-bialgebras \emph{(}Hopf $R$-algebras\emph{)}
with quasi-admissible filters $\frak{F}_{H},$ $\frak{F}_{K}$ and consider
the canonical $R$-linear mapping $\delta :H^{\ast }\otimes _{R}K^{\ast
}\rightarrow (H\otimes _{R}K)^{\ast }$ and the filter $\frak{F}$ of $R$%
-cofinite $H\otimes _{R}K$-ideals generated by $\frak{F}_{H}\times \frak{F}%
_{K}.$

\begin{enumerate}
\item  If $\frak{F}_{H}$ and $\frak{F}_{K}$ are moreover cofinitary (i.e.
admissible filters), then $(H\otimes _{R}K)_{\frak{F}}^{\circ }$ is an $R$%
-bialgebra \emph{(}respectively a Hopf $R$-algebra\emph{)}. If $R$ is
Noetherian, then $\delta $ induces an isomorphism of $R$-bialgebras \emph{(}%
respectively Hopf $R$-algebras\emph{) }$H_{\frak{F}_{H}}^{\circ }\otimes
_{R}K_{\frak{F}_{K}}^{\circ }\simeq (H\otimes _{R}K)_{\frak{F}}^{\circ }.$

\item  Let $R$ be Noetherian. If $\frak{F}_{K}$ is an $\alpha $-filter and $%
\frak{F}_{H}$ is cofinitary, then $(H\otimes _{R}K)_{\frak{F}}^{\circ }$ is
an $R$-bialgebra \emph{(}a Hopf $R$-algebra\emph{)} and $\delta $ induces an
isomorphism of $R$-bialgebras \emph{(}Hopf $R$-algebras\emph{) }$H_{\frak{F}%
_{H}}^{\circ }\otimes _{R}K_{\frak{F}_{K}}^{\circ }\simeq (H\otimes _{R}K)_{%
\frak{F}}^{\circ }.$
\end{enumerate}
\end{proposition}

\begin{definition}
The ring $R$ is called \emph{hereditary}, if every ideal $I\vartriangleleft R
$ is projective.
\end{definition}

\begin{theorem}
\label{erb-hopf}Let $R$ be Noetherian.

\begin{enumerate}
\item  If $H$ is an $\alpha $-bialgebra \emph{(}respectively a Hopf $\alpha $%
-algebra\emph{)}, then $(H,H^{\circ })\in \mathcal{P}_{Big}^{\alpha }$ \emph{%
(}respectively $(H,H^{\circ })\in \mathcal{P}_{Hopf}^{\alpha }$\emph{)}. If
moreover $H$ is commutative \emph{(}cocommutative\emph{)}, then $H^{\circ }$
is cocommutative \emph{(}commutative\emph{)}.

\item  If $R$ is hereditary, then there are \emph{self-adjoint}
contravariant functors
\begin{equation*}
\begin{tabular}{rrrrrrrrrrrr}
$(-)^{\circ }$ & $:$ & $\mathbf{Big}_{R}$ & $\rightarrow $ & $\mathbf{Big}%
_{R},$ & $(-)^{\circ }\text{ }$ & $:$ & $\mathbf{Hopf}_{R}$ & $\rightarrow $
& $\mathbf{Hopf}_{R}.$ &  &  \\
& $:$ & $\mathbf{CBig}_{R}$ & $\rightarrow $ & $\mathbf{CCBig}_{R},$ &  & $:$
& $\mathbf{CBig}_{R}$ & $\rightarrow $ & $\mathbf{CCBig}_{R}$ &  &  \\
& $:$ & $\mathbf{CCBig}_{R}$ & $\rightarrow $ & $\mathbf{CBig}_{R},$ &  & $:$
& $\mathbf{CCHopf}_{R}$ & $\rightarrow $ & $\mathbf{CHopf}_{R}$ &  &
\end{tabular}
\end{equation*}
\end{enumerate}
\end{theorem}

\begin{Beweis}
\begin{enumerate}
\item  If $H$ is an $\alpha $-bialgebra (a Hopf $\alpha $-algebra), then $%
H^{\circ }$ is by corollary \ref{zul-noeth} an $R$-bialgebra (a Hopf $R$%
-algebra) and $(H,H^{\circ })\in \mathcal{P}_{Big}^{\alpha }$ (respectively $%
(H,H^{\circ })\in \mathcal{P}_{Hopf}^{\alpha }$). The duality between the
commutativity and the cocommutativity follows now from \cite[Lemma 2.2]{Abu2}%
.

\item  Let $R$ be hereditary. Then for every $R$-bialgebra (respectively
Hopf $R$-algebra) $H$ the continuous dual $R$-module $H^{\circ }\subset R^{H}
$ is pure, \cite[Proposition 2.11]{AG-TW2000}, hence every $R$-bialgebra
(respectively Hopf $R$-algebra)\ is an $\alpha $-bialgebra (respectively a
Hopf $\alpha $-algebra) and $H^{\circ }$ is an $R$-bialgebra (respectively a
Hopf $R$-algebra). Moreover
\begin{equation*}
\Upsilon _{H,K}:\mathrm{Big}_{R}(H,K^{\circ })\rightarrow \mathrm{Big}%
_{R}(K,H^{\circ }),\text{ }f\mapsto \lbrack k\mapsto f(-)(k)]
\end{equation*}
is an isomorphism with inverse
\begin{equation*}
\Psi _{H,K}:\mathrm{Big}_{R}(K,H^{\circ })\rightarrow \mathrm{Big}%
_{R}(H,K^{\circ }),\text{ }g\mapsto \lbrack h\mapsto g(-)(h)].
\end{equation*}
It is easy to show that $\Upsilon _{H,K}$ and $\Psi _{H,K}$ are functorial
in $H$ and $K.\blacksquare $
\end{enumerate}
\end{Beweis}

\section{Coinduction functors in $\mathcal{P}_{Hopf}^{\protect\alpha }$}

In this section we consider the coinduction functors for the category of
Hopf $\alpha $-pairings respectively bialgebra $\alpha $-pairings that unify
several important situations (e.g. \cite{Don80}, \cite{APW91}, \cite[3.2]
{Lin93}).

\begin{definition}
Let $H$ be an $R$-bialgebra. For every left $H$-module $M$ we call the $R$%
-submodule
\begin{equation*}
M^{H}:=\{m\in M|\text{ }hm=\varepsilon (h)m\text{ for all }h\in H\}
\end{equation*}
the \emph{submodule of }$H$-\emph{invariants} of $M.$ For every right $H$%
-comodule $M$ we call
\begin{equation*}
M^{coH}:=\{m\in M|\text{ }\varrho _{M}(m)=m\otimes 1_{H}\}
\end{equation*}
the \emph{submodule of }$H$-\emph{coinvariants} of $M.$
\end{definition}

\begin{punto}
\label{big-md}Let $H$ be an $R$-bialgebra. If $M,N$ are right (respectively
left) $H$-modules, then $M\otimes _{R}N$ is a right (respectively a left) $H$%
-module with the \emph{canonical }$H$-module structure
\begin{equation}
(m\otimes n)h:=\sum mh_{1}\otimes nh_{2}\text{ (respectively }h(m\otimes
n):=\sum h_{1}m\otimes h_{2}n\text{).}  \label{right-MN}
\end{equation}
In particular the ground ring $R$ is an $H$-bimodule through
\begin{equation*}
h\rightharpoonup r:=\varepsilon (h)r=:r\leftharpoonup h\text{ for all }h\in H%
\text{ and }r\in R.
\end{equation*}
\end{punto}

\begin{punto}
\label{Big-com}Let $K$ be an $R$-bialgebra. If $M,N$ are right (respectively
left) $K$-comodules, then $M\otimes _{R}N$ is a right (respectively a left) $%
K$-comodule through the \emph{canonical }right (respectively left) $K$%
-comodule structure
\begin{equation}
m\otimes n\mapsto \sum m_{<0>}\otimes n_{<0>}\otimes m_{<1>}n_{<1>}\text{\
(resp. }m\otimes n\mapsto \sum m_{<-1>}n_{<-1>}\otimes m_{<0>}n_{<0>}\text{).%
}  \label{bial-com}
\end{equation}
In particular the ground ring $R$ is a $K$-bicomodule throughout
\begin{equation*}
R\rightarrow R\otimes _{R}K,\text{ }r\mapsto r\otimes 1_{K}\text{ and }%
R\rightarrow K\otimes _{R}R,\text{ }r\mapsto 1_{K}\otimes r.
\end{equation*}
\end{punto}

\begin{lemma}
\label{coK=H}Let $P=(H,K)\in \mathcal{P}_{Big},$ $(M,\varrho _{M})$ be a
right $K$-comodule and consider $M$ with the induced left $H$-module
structure. If $\alpha _{M}^{P}:M\otimes _{R}K\rightarrow \mathrm{Hom}%
_{R}(H,M)$ is injective, then $M^{H}=M^{coK}.$
\end{lemma}

\begin{Beweis}
We have for all $m\in M^{coK}$ and $h\in H:$%
\begin{equation*}
hm=m<h,1_{K}>=m\varepsilon _{H}(h)\text{ for every }h\in H,
\end{equation*}
i.e. $m\in M^{H}.$ On the other hand, we have for all $m\in M^{H}$ and $h\in
H:$%
\begin{equation*}
\begin{tabular}{lllll}
$\alpha _{M}^{P}(\sum m_{<0>}\otimes m_{<1>})(h)$ & $=$ & $\sum
m_{<0>}<h,m_{<1>}>$ & $=$ & $hm$ \\
& $=$ & $m\varepsilon _{H}(h)$ & $=$ & $m<h,1_{K}>$ \\
& $=$ & $\alpha _{M}^{P}(m\otimes 1_{K})(h).$ &  &
\end{tabular}
\end{equation*}
If $\alpha _{M}^{P}$ is injective, then $\varrho _{M}(m)=\sum m_{<0>}\otimes
m_{<1>}=m\otimes 1_{K},$ i.e. $m\in M^{coK}$ and consequently $%
M^{H}=M^{coK}.\blacksquare $
\end{Beweis}

\begin{lemma}
\label{_H=H}Let $H$ be a Hopf $R$-algebra and $M,N\in $ $_{H}\mathcal{M}.$
Then $\mathrm{Hom}_{R}(M,N)$ is a left $H$-module through
\begin{equation}
(hf)(m)=\sum h_{1}f(S_{H}(h_{2})m)\text{ for all }h\in H,m\in M\text{ and }%
f\in \mathrm{Hom}_{R}(M,N).  \label{ind-mod}
\end{equation}
Moreover $\mathrm{Hom}_{H-}(M,N)=\mathrm{Hom}_{R}(M,N)^{H}.$
\end{lemma}

\begin{Beweis}
For all $h,\widetilde{h}\in H,$ $f\in \mathrm{Hom}_{R}(M,N)$ and $m\in M$ we
have
\begin{equation*}
\begin{tabular}{lllll}
$((h\widetilde{h})f)(m)$ & $:=$ & $\sum (h\widetilde{h})_{1}f(S_{H}((h%
\widetilde{h})_{2})m)$ & $=$ & $\sum h_{1}\widetilde{h}_{1}f(S_{H}(h_{2}%
\widetilde{h}_{2})m)$ \\
& $=$ & $\sum h_{1}\widetilde{h}_{1}f(S_{H}(\widetilde{h}_{2})S_{H}(h_{2})m)$
& $=$ & $\sum h_{1}((\widetilde{h}f)(S_{H}(h_{2})m))$ \\
& $=$ & $(h(\widetilde{h}f))(m),$ &  &
\end{tabular}
\end{equation*}
i.e. $\mathrm{Hom}_{R}(M,N)$ is a left $H$-module with the left $H$-action (%
\ref{ind-mod}).

For all $f\in \mathrm{Hom}_{H-}(M,N),$ $h\in H$ and $m\in M$ we have
\begin{equation*}
(hf)(m):=\sum h_{1}f(S_{H}(h_{2})m)=\sum h_{1}S_{H}(h_{2})f(m)=(\varepsilon
(h)1_{H})f(m)=(\varepsilon (h)f)(m),
\end{equation*}
i.e. $f\in \mathrm{Hom}_{R}(M,N)^{H}.$ On the other hand, if $g\in \mathrm{%
Hom}_{R}(M,N)^{H},$ then we have for all $h\in H$ and $m\in M:$%
\begin{equation*}
\begin{tabular}{lllll}
$g(hm)$ & $=$ & $g(\sum \varepsilon (h_{1})h_{2}m)$ & $=$ & $\sum
(\varepsilon (h_{1})g)(h_{2}m)$ \\
& $=$ & $\sum (h_{1}g)(h_{2}m)$ & $=$ & $\sum h_{11}(g(S_{H}(h_{12})h_{2}m))$
\\
& $=$ & $\sum h_{1}(g(S_{H}(h_{21})h_{22}m))$ & $=$ & $\sum
h_{1}g(\varepsilon (h_{2})1_{H}m)$ \\
& $=$ & $(\sum h_{1}\varepsilon (h_{2}))g(m)$ & $=$ & $hg(m),$%
\end{tabular}
\end{equation*}
i.e. $g\in \mathrm{Hom}_{H-}(M,N).\blacksquare $
\end{Beweis}

\qquad The following lemma generalizes the corresponding results \cite[Page
165]{Lin93} and \cite[Page 103]{Lin94}:

\begin{lemma}
\label{H=coK}Let $P=(H,K),$ $Q=(Y,Z)\in \mathcal{P}_{Hopf}^{\alpha },$ $(\xi
,\theta ):(Y,Z)\rightarrow (H,K)$ be a morphism in $\mathcal{P}%
_{Hopf}^{\alpha }$ and $N\in $ $_{H}\mathcal{M}.$

\begin{enumerate}
\item  $\mathrm{Hom}_{R}(Y,N)$ is a left $H$-module through
\begin{equation}
(hf)(y)=\sum h_{1}f(S_{Y}(\xi (h_{2})y))\text{ for all }h\in H,f\in \mathrm{%
Hom}_{R}(Y,N)\text{ and }y\in Y.  \label{(Y,N)}
\end{equation}

\item  If we consider $\mathrm{Hom}_{R}(Y,N)$ with the canonical left $Y$%
-module structure, then
\begin{equation*}
h({\normalsize y}f)={\normalsize y}(hf)\text{ for all }h\in H,{\normalsize y}%
\in {\normalsize Y}\text{ and }f\in \mathrm{Hom}_{R}({\normalsize Y},N).
\end{equation*}
So $\mathrm{Hom}_{R}(Y,N)^{H}\subseteq \mathrm{Hom}_{R}(Y,N)$ is a left $Y$%
-submodule.

\item  If $_{H}N$ is $K$-rational, then $N\otimes _{R}Z$ is a right $K$%
-comodule through
\begin{equation}
\psi :N\otimes _{R}Z\rightarrow N\otimes _{R}Z\otimes _{R}K,\text{ }n\otimes
z\mapsto \sum n_{<0>}\otimes z_{2}\otimes n_{<1>}S_{K}(\theta (z_{1})).
\label{NZ}
\end{equation}
\end{enumerate}
\end{lemma}

\begin{Beweis}
\begin{enumerate}
\item  By assumption $\xi :H\rightarrow Y$ is a Hopf $R$-algebra morphism
and so $\xi (S_{H}(h))=S_{Y}(\xi (h))$ for every $h\in H.$ If we consider
the left $H$-module $Y_{\xi },$ then the left $H$-action on $\mathrm{Hom}%
_{R}(Y_{\xi },N)$ in (\ref{ind-mod}) coincides with that in (\ref{(Y,N)}),
hence $\mathrm{Hom}_{R}(Y,N)$ is a left $H$-module by Lemma \ref{_H=H}.

\item  Trivial.

\item  $Z$ is obviously a right $K$-comodule through
\begin{equation*}
\varrho _{Z}:Z\rightarrow Z\otimes _{R}K,\text{ }z\mapsto \sum z_{2}\otimes
S_{K}(\theta (z_{1}))\text{ for every }z\in Z.
\end{equation*}
By assumption and Theorem \ref{equal} $N$ is a right $K$-comodule and so $%
(N\otimes _{R}Z,\psi )$ is, by \ref{Big-com}, a right $K$-comodule.$%
\blacksquare $
\end{enumerate}
\end{Beweis}

\begin{punto}
\label{Coind-Hopf}Let $P=(H,K),$ $Q=(Z,Y)\in \mathcal{P}_{Hopf}^{\alpha }$
and $(\xi ,\theta ):(Y,Z)\rightarrow (H,K)$ be a morphism in $\mathcal{P}%
_{Hopf}^{\alpha }.$ For every $N\in \mathrm{Rat}^{K}(_{H}\mathcal{M})$
consider $N\otimes _{R}Z$ with the right $K$-comodule structure (\ref{NZ}).
If we consider the coinduction functor
\begin{equation*}
\mathrm{Coind}_{P}^{Q}(-):\mathrm{Rat}{\normalsize ^{K}(_{H}\mathcal{M}%
)\rightarrow \mathrm{Rat}^{Z}(_{Y}}\mathcal{M}{\normalsize ),}\text{ }%
N\mapsto \mathrm{HOM}_{H-}(Y,N):={\normalsize \mathrm{Rat}^{Z}(_{Y}(}\mathrm{%
Hom}_{H-}(Y,N))),
\end{equation*}
then we have functorial isomorphisms
\begin{equation*}
\begin{tabular}{lllll}
$(N\otimes _{R}Z)^{coK}$ & $\simeq $ & $(N\otimes _{R}Z)^{H}$ &  & (Lemma
\ref{coK=H}); \\
& $\simeq $ & $\mathrm{HOM}_{R}(Y,N)^{H}$ &  & (\ref{otC=rtc}); \\
& $=$ & $\mathrm{Rat}^{Z}(_{Y}(\mathrm{Hom}_{R}(Y,N)^{H}))$ &  &  \\
& $=$ & $\mathrm{HOM}_{H-}(Y,N):=\mathrm{Coind}_{P}^{Q}(N)$ &  & (Lemma \ref
{_H=H}); \\
& $\simeq $ & $N\square _{K}Z$ &  & (Theorem \ref{phi-ad}); \\
& $\simeq $ & $\mathrm{Hom}_{H^{e}}(H,N\otimes _{R}Z).$ &  & (Proposition
\ref{cot=Hom(ot)}).
\end{tabular}
\end{equation*}
\end{punto}

\begin{corollary}
\label{Coind-ot}Let $P=(H,K),$ $Q=(Y,Z)\in \mathcal{P}_{Hopf}^{\alpha }$ and
$(\xi ,\theta ):(Y,Z)\rightarrow (H,K)$ be a morphism in $\mathcal{P}%
_{Hopf}^{\alpha }.$ Let $M\in $ $\mathcal{M}^{Z},$ $N\in \mathcal{M}^{K}$
and consider $M^{\theta }\otimes _{R}N$ with the canonical right $K$%
-comodule structure. If $M_{R}$ is flat, then there is an isomorphism of $Z$%
-comodules
\begin{equation*}
\mathrm{Coind}_{P}^{{\normalsize Q}}(M^{\theta }\otimes _{R}N)\simeq
(M^{\theta }\otimes _{R}N)\square _{K}{\normalsize Z}\simeq M\otimes
_{R}(N\square _{K}{\normalsize Z})\simeq M\otimes _{R}\mathrm{Coind}_{P}^{%
{\normalsize Q}}(N).
\end{equation*}
\end{corollary}

\section{Classical Duality}

\qquad Over a commutative base field one has a \emph{duality} between the
groups and the commutative Hopf algebras (e.g. \cite[9.3]{Mon93}, \cite
{Swe69}). In this section we show that such a duality is valid over \emph{%
hereditary\ Noetherian }ground rings.

\begin{definition}
Let $(C,\Delta ,\varepsilon )$ be an $R$-coalgebra. With
\begin{equation*}
\mathcal{G}(C):=\{0\neq x\in C|\text{ }\Delta (x)=x\otimes x\text{ and }%
\varepsilon (x)=1_{R}\}
\end{equation*}
we denote the set of \emph{group-like elements} of $C.$ If $x,y\in \mathcal{G%
}(C),$ then we denote with
\begin{equation*}
P_{(x,y)}(c):=\{c\in C|\text{ }\Delta (c)=x\otimes c+c\otimes y\}
\end{equation*}
the set of $(x,y)$\emph{-primitive elements} in $C.$ For an $R$-bialgebra $B$
we call the $(1_{B},1_{B})$-primitive elements of $B$ \emph{primitive
elements.}
\end{definition}

\qquad The following result is easy to prove

\begin{lemma}
\label{gr-like}Let $C$ be an $R$-coalgebra.

\begin{enumerate}
\item  If $D$ is an $R$-coalgebra and $f:D\rightarrow C$ is an $R$-coalgebra
morphism, then $f(\mathcal{G}(D))\subseteq \mathcal{G}(C).$

\item  If $\{0_{R},1_{R}\}$ are the only idempotents in $R$ \emph{(}e.g. $R$
is a domain\emph{)} and $\Delta _{C}(x)=x\otimes x$ for some $0\neq x\in C,$
then $\varepsilon _{C}(x)=1_{R},$ i.e. $x\in \mathcal{G}(C).$

\item  If $x,y\in \mathcal{G}(C)$ and $c\in P_{(x,y)}(C),$ then $\varepsilon
_{C}(c)=0.$

\item  For every $R$-coalgebra $C$ we have a bijection
\begin{equation*}
\mathrm{Cog}_{R}(R,C)\leftrightarrow \mathcal{G}(C),\text{ }f\mapsto f(1_{R})%
\text{ and }x\mapsto \lbrack 1_{R}\mapsto x]\text{ }\forall \text{ }f\in
\mathrm{Cog}_{R}(R,C),\text{ }x\in \mathcal{G}(C).
\end{equation*}

\item  If $R$ is Noetherian and $A$ is an $\alpha $-algebra, then $\mathrm{%
Alg}_{R}(A,R)=\mathcal{G}(A^{\circ })=\mathrm{Cog}_{R}(R,A^{\circ }).$
\end{enumerate}
\end{lemma}

\begin{punto}
\label{kokomm}For every set $G$ the free $R$-module $RG$ becomes a \emph{%
cocommutative }$R$-coalgebra $\mathcal{K}(G):=(RG,\Delta _{g},\varepsilon
_{g}),$ where the comultiplication $\Delta _{g}$ and the counit $\varepsilon
_{g}$ are given by the linear extension of their images on the elements of $%
G:$%
\begin{equation*}
\Delta _{g}(x)=x\otimes x\text{ and }\varepsilon _{g}(x)=1\text{ for every }%
x\in G.
\end{equation*}
If $(G,\mu _{G},e_{G})$ is a monoid, then $\mu _{G}$ respectively $e_{G}$
induce on $RG$ a multiplication $\mu $ respectively a unity $\eta ,$ such
that $\mathcal{K}(G)=(RG,\mu ,\eta ,\Delta _{g},\varepsilon _{g})$ is an $R$%
-bialgebra. If $G$ is moreover a group, then $RG$ is a Hopf $R$-algebra with
antipode defined on the basis elements as $S_{g}:RG\rightarrow RG,$ $%
x\mapsto x^{-1}$ for every $x\in G.$ On the other hand, let $H$ be an $R$%
-bialgebra. Then $\Delta _{H}(1_{H})=1_{H}\otimes 1_{H}$ and we have for all
$x,y\in \mathcal{G}(H):$%
\begin{equation*}
\Delta _{H}(xy)=\Delta _{H}(x)\Delta _{H}(y)=(x\otimes x)(y\otimes
y)=xy\otimes xy,
\end{equation*}
i.e. $xy$ is a group-like element in $H$ and $\mathcal{G}(H)$ is a monoid.
If $H$ is moreover a Hopf $R$-algebra and $x\in \mathcal{G}(H),$ then $%
x^{-1}:=S_{H}(x)\in \mathcal{G}(H),$ i.e. $\mathcal{G}(H)$ is a group.
\end{punto}

\begin{proposition}
\label{koko-ad}\emph{(\cite{Gru69})} Denote with $\mathbf{Ens},$ $\mathbf{Mon%
}$ and $\mathbf{Gr}$ the categories of sets, monoids and groups
respectively. Then we have adjoint pairs of covariant functors $(\mathcal{K}%
(-),\mathcal{G}(-)):$%
\begin{equation*}
\begin{tabular}{lllllllllll}
$\mathcal{K}(-)$ & $:$ & $\mathbf{Ens}$ & $\rightarrow $ & $\mathbf{CCog}%
_{R},\ $ & $\mathcal{G}(-)$ & $:$ & $\mathbf{CCog}_{R}$ & $\rightarrow $ & $%
\mathbf{Ens}$ &  \\
& $:$ & $\mathbf{Mon}$ & $\rightarrow $ & $\mathbf{CCBialg}_{R},$ &  & $:$ &
$\mathbf{CCBialg}_{R}$ & $\rightarrow $ & $\mathbf{Mon}$ &  \\
& $:$ & $\mathbf{Gr}$ & $\rightarrow $ & $\mathbf{CCHopf}_{R},$ &  & $:$ & $%
\mathbf{CCHopf}_{R}$ & $\rightarrow $ & $\mathbf{Gr.}$ &
\end{tabular}
\end{equation*}
If $R$ is moreover an integral domain, then we have a natural isomorphism $%
\mathcal{G}(-)\circ \mathcal{K}(-)\simeq id.$
\end{proposition}

\subsection*{Representative mappings}

\begin{punto}
Let $R$ be Noetherian, $(G,\mu ,e)$ be a monoid (respectively a group) and
denote with
\begin{equation*}
\mathcal{R}(G):=\{f\in R^{G}|\text{ }GfG\text{ is finitely generated in }%
\mathcal{M}_{R}{\normalsize \}\simeq }(RG)^{\circ }
\end{equation*}
the set of \emph{representative mappings }on $G.$ We call $G$ an $\alpha $%
\emph{-monoid }(respectively an $\alpha $\emph{-group}), if $(RG,\mathcal{R}%
(G))$ is an $\alpha $-pairing, or equivalently if $\mathcal{R}(G)\subset
R^{G}$ is pure.
\end{punto}

As a consequence of Lemma \ref{A0=} and Corollary \ref{zul-noeth} we get

\begin{corollary}
\label{repres}Let $R$ be Noetherian. If $G$ is an $\alpha $-monoid, then $%
\mathcal{R}(G)$ is an $R$-bialgebra. If $G$ is moreover an $\alpha $-group,
then $\mathcal{R}(G)$ is a Hopf $R$-algebra with antipode
\begin{equation*}
S:\mathcal{R}{\normalsize (G)\rightarrow }\mathcal{R}{\normalsize (G),}\text{
}S(f)(x)=f(x^{-1})\text{ for }f\in \mathcal{R}{\normalsize (G)}\text{ and }%
{\normalsize x\in G.}
\end{equation*}
\end{corollary}

\begin{Notation}
Let $G$ be a monoid. The category of unital left (respectively right) $G$%
-modules is denoted by $_{G}\mathcal{M}$ (respectively $\mathcal{M}_{G}$).
\end{Notation}

As a consequence of Theorem \ref{R(G)-co} we get

\begin{corollary}
\label{C_R(G)}Let $R$ be Noetherian, $G$ be a monoid and $C\subseteq
\mathcal{R}(G)$ be a $G$-subbimodule. If $P=(RG,C)$ is an $\alpha $-pairing,
then $C$ is an $R$-coalgebra and we have category isomorphisms
\begin{equation*}
\begin{tabular}{lllll}
$\mathcal{M}^{C}$ & $\simeq $ & $\mathrm{Rat}^{C}(_{G}\mathcal{M})$ & $=$ & $%
\sigma \lbrack _{RG}C]$ \\
& $\simeq $ & $\mathrm{Rat}^{C}(_{C^{\ast }}\mathcal{M})$ & $=$ & $\sigma
\lbrack _{C^{\ast }}C]$%
\end{tabular}
\&
\begin{tabular}{lllll}
$^{C}\mathcal{M}$ & $\simeq $ & $^{C}\mathrm{Rat}(\mathcal{M}_{G})$ & $=$ & $%
\sigma \lbrack C_{RG}]$ \\
& $\simeq $ & $^{C}\mathrm{Rat}(\mathcal{M}_{C^{\ast }})$ & $=$ & $\sigma
\lbrack C_{C^{\ast }}].$%
\end{tabular}
\end{equation*}
\end{corollary}

\begin{punto}
Let $G$ be a monoid. A left (respectively right) $G$-module will be called
\emph{locally finite}, if $(RG)m$ (respectively $m(RG)$) is finitely
generated in $\mathcal{M}_{R}$ for every $m\in M.$ For every monoid $G$
denote with $\mathrm{Loc}(_{G}\mathcal{M})\subset $ $_{G}\mathcal{M}$
(respectively $\mathrm{Loc}(\mathcal{M}_{G})\subseteq \mathcal{M}_{G}$) the
full subcategory of locally finite left (respectively right) $G$-modules.
\end{punto}

\qquad As a consequence of \cite[Proposition 3.23]{Abu2} we get

\begin{proposition}
Let $R$ be Noetherian and $G$ be a monoid.

\begin{enumerate}
\item  Every $\mathcal{R}(G)$-subgenerated left \emph{(}respectively right%
\emph{)} $G$-module is locally finite.

\item  If $RG$ is cofinitely $R$-cogenerated, then $\sigma \lbrack _{G}%
\mathcal{R}(G)]=\mathrm{Loc}(_{G}\mathcal{M})$ and $\sigma \lbrack \mathcal{R%
}(G)_{G}]=\mathrm{Loc}(\mathcal{M}_{G}).$ If $G$ is moreover an $\alpha $%
-monoid, then we have category isomorphisms
\begin{equation*}
\begin{tabular}{lllllll}
$\mathcal{M}^{\mathcal{R}(G)}$ & $\simeq $ & $\mathrm{Rat}^{\mathcal{R}%
(G)}(_{G}\mathcal{M})$ & $=$ & $\sigma \lbrack _{G}\mathcal{R}(G)]$ & $=$ & $%
\mathrm{Loc}(_{G}\mathcal{M});$ \\
$^{\mathcal{R}(G)}\mathcal{M}$ & $\simeq $ & $^{\mathcal{R}(G)}\mathrm{Rat}(%
\mathcal{M}_{G})$ & $=$ & $\sigma \lbrack \mathcal{R}(G)_{G}]$ & $=$ & $%
\mathrm{Loc}(\mathcal{M}_{G}).$%
\end{tabular}
\end{equation*}
\end{enumerate}
\end{proposition}

\qquad The following result generalizes the classical duality between
monoids (groups) and \emph{commutative} $R$-bialgebras (Hopf $R$-algebras),
e.g. \cite[9.3]{Mon93}, from the case of base fields to the case of
arbitrary hereditary Noetherian rings.

\begin{theorem}
\label{ad-hopf-gr}If $R$ is Noetherian and hereditary, then there is a \emph{%
duality} between monoids \emph{(}respectively groups\emph{)} and commutative
$R$-bialgebras \emph{(}respectively Hopf $R$-algebras\emph{)} through the
right-adjoint contravariant functors
\begin{equation*}
\begin{tabular}{llllllllll}
$\mathcal{R}(-)$ & $:$ & $\mathbf{Mon}$ & $\rightarrow $ & $\mathbf{CBig}%
_{R},$ & $\mathrm{Alg}_{R}(-,R)$ & $:$ & $\mathbf{CBig}_{R}$ & $\rightarrow $
& $\mathbf{Mon},$ \\
& $:$ & $\mathbf{Gr}$ & $\rightarrow $ & $\mathbf{CHopf}_{R},$ &  & $:$ & $%
\mathbf{CHopf}_{R}$ & $\rightarrow $ & $\mathbf{Gr}.$%
\end{tabular}
\end{equation*}
\end{theorem}

\begin{Beweis}
Let $R$ be Noetherian and hereditary. Then for every $R$-algebra $A,$ the
character module $A^{\circ }\subset R^{A}$ is pure (e.g. \cite[Proposition
2.11]{AG-TW2000}). If $G$ is a monoid (respectively a group), then $\mathcal{%
K}(G)=(RG,\mu ,\eta ,\Delta _{g},\varepsilon _{g})$ is by \ref{kokomm} a
\emph{cocommutative}{\normalsize \ }$R$-bialgebra (respectively Hopf $R$%
-algebra) and so $\mathcal{R}(G)=(RG)^{\circ }$ is by Theorem \ref{erb-hopf}
a \emph{commutative }$R$-bialgebra (respectively Hopf $R$-algebra). If $H$
is an $R$-bialgebra (respectively a Hopf $R$-algebra), then $H^{\circ }$ is
by Theorem \ref{erb-hopf} an $R$-bialgebra (respectively a Hopf $R$%
-algebra), hence $\mathrm{Alg}_{R}(H,R)=\mathcal{G}(H^{\circ })\ $is a
monoid (respectively a group). It is easy to see that we have isomorphisms
of functors
\begin{equation*}
\mathcal{R}{\normalsize (-)\simeq (-)^{\circ }\circ \mathcal{K}(-)}\text{
and }\mathrm{Alg}_{R}(-,R)\simeq \mathcal{G}{\normalsize (-)\circ (-)^{\circ
}.}
\end{equation*}
The result follows now from Theorems \ref{erb-hopf} and \ref{koko-ad}.$%
\blacksquare $
\end{Beweis}

\section{Affine group schemes}

\emph{Affine groups schemes} over arbitrary commutative ground rings were
presented by J. Jantzen \cite{Jan87}. If $\frak{G}$ is an affine group
scheme with \emph{coordinate ring }$R(\frak{G})$, then the category of \emph{%
left }$\frak{G}$\emph{-modules} $_{\frak{G}}\mathcal{M}$ and the category of
right $R(\frak{G})$-comodules $\mathcal{M}^{R(\frak{G})}$ are equivalent. In
the case $R(\frak{G})$ is locally projective as an $R$-module we extend this
equivalence to the category of $R(\frak{G})$-rational left $R(\frak{G}%
)^{\ast }$-modules $\mathrm{Rat}^{R(\frak{G})}(_{R(\frak{G})^{\ast }}%
\mathcal{M})$ which turns to be equal to the category of $R(\frak{G})$%
-subgenerated left $R(\frak{G})^{\ast }$-modules $\sigma \lbrack _{R(\frak{G}%
)^{\ast }}R(\frak{G})].$ It follows that in this case $_{\frak{G}}\mathcal{M}
$ is a Grothendieck category of type $\sigma \lbrack M]$ and one can use the
well developed theory of such categories (e.g. \cite{Wis88}, \cite{Wis96})
to study the category $_{\frak{G}}\mathcal{M}.$

\begin{punto}
With an $R$\emph{-functor} (respectively\textbf{\ }a \emph{monoid }$R$\emph{%
-functor, }a \emph{group }$R$\emph{-functor}) we understand a functor from
the category of \emph{commutative }$R$-algebras $\mathbf{CAlg}_{R}$ to $%
\mathbf{Ens}$ (respectively to $\mathbf{Mon},$ $\mathbf{Gr}$). An \emph{%
affine scheme}\textbf{\ }(respectively an \emph{affine monoid scheme},%
{\normalsize \ }an \emph{affine group scheme}\textbf{)} over $R$ is a \emph{%
representable} $R$-functor (respectively monoid $R$-functor, group $R$%
-functor)
\begin{equation*}
\begin{tabular}{lllll}
$\frak{G}=\mathrm{Alg}_{R}(H,-)$ & $:$ & $\mathbf{CAlg}_{R}$ & $\rightarrow $
& $\mathbf{Ens},$ \\
& $:$ & $\mathbf{CBig}_{R}$ & $\rightarrow $ & $\mathbf{Mon},$ \\
& $:$ & $\mathbf{CHopf_{R}}$ & $\rightarrow $ & $\mathbf{Gr}.$%
\end{tabular}
\end{equation*}
The \emph{commutative} $R$-algebra $H$ is called the \emph{coordinate ring}
of $\frak{G}$ and is denoted with $R(\frak{G}).$ With $\mathbf{Aff}_{R}$
(respectively $\mathbf{AffMon}_{R},$ $\mathbf{AffGr}_{R}$) we denote the
category of affine schemes (respectively affine monoid schemes, affine group
schemes) with morphisms the \emph{natural transformations.}
\end{punto}

\begin{punto}
$\frak{G}$\textbf{-modules.} (\cite[2.7]{Jan87}) Let $\frak{G}=\mathrm{Alg}%
_{R}(H,-)$ be an affine group scheme. An $R$-module $M$ is said to be a
\emph{left }(respectively \emph{a right}) $\frak{G}$\emph{-module}, if there
is a $\frak{G}(A)$ module structure on $M\otimes _{R}A$ (respectively on $%
A\otimes _{R}M$), functorial in $A,$ for every commutative $R$-algebra $A.$
The category of left (respectively right) $\frak{G}$-modules and $\frak{G}$%
-linear mappings will be denoted by $_{\frak{G}}\mathcal{M}$ (respectively
by $\mathcal{M}_{\frak{G}}$).
\end{punto}

\begin{punto}
\textbf{Yoneda Lemma.}\label{Yoneda} (\cite[44.3]{Wis88}) Let $\frak{C}$ be
a category, $F:\frak{C}\rightarrow \mathbf{Ens}$ be a covariant functor and
denote for $A\in \frak{C}$ the class of functorial morphisms between $Mor_{%
\frak{C}}(A,-)$ and $F$ with $\underline{\underline{Nat}}(Mor_{\frak{C}%
}(A,-),F).$ Then the following \emph{Yoneda-mapping} is bijective:
\begin{equation*}
\underline{\underline{Nat}}(Mor_{\frak{C}}(A,-),F)\rightarrow F(A),\text{ }%
\phi \mapsto \phi _{A}(id_{A}).
\end{equation*}
\end{punto}

With the help of Yoneda-Lemma (Compare \cite[Chapter 2]{Jan87}) one obtains:

\begin{proposition}
\label{affGr=cHop} Let $R$ be an arbitrary commutative ring.

\begin{enumerate}
\item  If $\frak{G}=\mathrm{Alg}_{R}(H,-)$ is an affine monoid scheme \emph{(%
}respectively an affine group scheme\emph{)}, then the coordinate ring $H=R(%
\frak{G})$ is an $R$-bialgebra \emph{(}respectively a Hopf $R$-algerba\emph{)%
} and we have equivalences of categories
\begin{equation*}
\mathbf{AffMon}_{R}\approx (\mathbf{CBig}_{R})^{op}\text{ and }\mathbf{AffGr}%
_{R}\approx (\mathbf{CHopf}_{R})^{op}.
\end{equation*}

\item  For every affine group scheme $\frak{G}$ with coordinate ring $R(%
\frak{G}),$ the category of left $\frak{G}$-modules $_{\frak{G}}\mathcal{M}$
and the category of right $R(\frak{G})$-comodules $\mathcal{M}^{R(\frak{G})}$
are equivalent.
\end{enumerate}
\end{proposition}

\begin{punto}
Let $\frak{G}$ be an affine group scheme with coordinate ring $R(\frak{G}),$
$\omega :=\mathrm{Ker}(\varepsilon _{R(\frak{G})}),$ $\frak{F}_{\omega
}:=\{\omega ^{n}|$ $n\geq 1\}$ and consider $R(\frak{G})^{\ast }$ with the
finite topology and $R(\frak{G})$ with the induced left linear topology $%
\frak{T}(\frak{F}_{\omega }).$ By \cite[7.7]{Jan87}
\begin{equation}
hy(\frak{G}):=\{f\in R(\frak{G})^{\ast }|\text{ }f(\omega ^{n})=0\text{ for
some }n\geq 1\}  \label{hy(G)}
\end{equation}
is an $R$-subalgebra of $R(\frak{G})^{\ast },$ the so called \emph{%
hyperalgebra} of $\frak{G},$ and we get a measuring $R$-pairing $(hy(\frak{G}%
),R(\frak{G})).$ If $hy(\frak{G})\subset R(\frak{G})^{\ast }$ is dense, then
we call $\frak{G}$ \emph{connected}. If $R(\frak{G})/\omega ^{n}$ is
finitely generated projective in $\mathcal{M}_{R}$ for every $n\geq 1,$ then
$\frak{G}$ is called \emph{infinitesimal flat}. We say $\frak{G}$ satisfies
the $\alpha $\emph{-condition} (or $\frak{G}$ is an affine $\alpha $\emph{%
-group scheme}), if $(hy(\frak{G}),R(\frak{G}))$ satisfies the $\alpha $%
-condition. We call $\frak{G}\;$\emph{locally projective}, if $R(\frak{G})$
is locally projective as an $R$-module.
\end{punto}

\begin{theorem}
\label{Gmod=}Let $\frak{G}$ be an affine group scheme with coordinate ring $%
R(\frak{G}).$

\begin{enumerate}
\item  If $\frak{G}$ is locally projective, then there are equivalences of
categories
\begin{equation*}
_{\frak{G}}\mathcal{M}\approx \mathcal{M}^{R(\frak{G})}\simeq \mathrm{Rat}%
^{R(\frak{G})}(_{R(\frak{G})^{\ast }}\mathcal{M})=\sigma \lbrack _{R(\frak{G}%
)^{\ast }}R(\frak{G})].
\end{equation*}

\item  $\frak{G}$ is an affine $\alpha $-group scheme if and only if $\frak{G%
}$ is locally projective and connected. If these equivalent conditions are
satisfied, then we have equivalences of categories
\begin{equation*}
\begin{tabular}{lllllll}
$_{\frak{G}}\mathcal{M}$ & $\approx $ & $\mathcal{M}^{R(\frak{G})}$ & $%
\simeq $ & $\mathrm{Rat}^{R(\frak{G})}(_{R(\frak{G})^{\ast }}\mathcal{M})$ &
$=$ & $\sigma \lbrack _{R(\frak{G})^{\ast }}R(\frak{G})]$ \\
&  &  & $\simeq $ & $\mathrm{Rat}^{R(\frak{G})}(_{hy(\frak{G})}\mathcal{M})$
& $=$ & $\sigma \lbrack _{hy(\frak{G})}R(\frak{G})].$%
\end{tabular}
\end{equation*}

\item  The following are equivalent:

(i) $\frak{G}$ is connected \emph{(}i.e. $hy(\frak{G})\subset R(\frak{G}%
)^{\ast }$ is dense\emph{)};

(ii) $\sigma \lbrack _{hy(\frak{G})}R(\frak{G})]=\sigma \lbrack _{R(\frak{G}%
)^{\ast }}R(\frak{G})].$

If $R$ is a injective cogenerator, then (i), (ii) are moreover equivalent to:

(iii) $R(\frak{G})\hookrightarrow hy(\frak{G})^{\ast };$

(iv) $\frak{T}(\frak{F}_{\omega })$ is Hausdorff.
\end{enumerate}
\end{theorem}

\begin{Beweis}
\begin{enumerate}
\item  The equivalence $_{\frak{G}}\mathcal{M}\approx \mathcal{M}^{R(\frak{G}%
)}$ follows from Proposition \ref{affGr=cHop}. The remaining category
isomorphisms follow from Theorem \ref{equal}.

\item  Follows from Theorem \ref{equal}.

\item  $hy(\frak{G})\subset R(\frak{G})^{\ast }$ is an $R$-subalgebra and so
the equivalence \emph{(i)} $\Leftrightarrow $ \emph{(ii)} follows by Lemma
\ref{dense}.

Let $R$ be an injective cogenerator.

The equivalence \emph{(i)} $\Leftrightarrow $ \emph{(iii)} follows from
\cite[Theorem 1.8 (2)]{Abu1}. Consider now the measuring $R$-pairings $\frak{%
G}:=(hy(\frak{G}),R(\frak{G})).$ Then we have
\begin{equation*}
\overline{0_{R(\frak{G})}}=\bigcap_{n=1}^{\infty }\omega
^{n}=\bigcap_{n=1}^{\infty }\mathrm{KeAn}(\omega ^{n})=\mathrm{Ke}%
(\sum_{n=1}^{\infty }\mathrm{An}(\omega ^{n}))=\mathrm{Ke}(hy(\frak{G}))=%
\mathrm{Ker}(\chi _{\frak{G}}).
\end{equation*}
Consequently $\frak{T}(\frak{F}_{\omega })$ is Hausdorff if and only if $R(%
\frak{G})\overset{\chi _{\frak{G}}}{\hookrightarrow }hy(\frak{G})^{\ast }$
and we are done.$\blacksquare $
\end{enumerate}
\end{Beweis}

\section*{Coinduction functors for affine $\protect\alpha $-schemes}

\begin{punto}
Let $\frak{G},$ $\frak{H}$ be affine $\alpha $-group schemes and $\varphi :%
\frak{H}\rightarrow \frak{G}$ be a morphism in $\mathbf{AffGr}_{R}.$ Then $%
\varphi $ induces a Hopf $R$-algebra morphism $\varphi _{\#}:R(\frak{G}%
)\rightarrow $\emph{\ }$R(\frak{H})$ (called a \emph{comorphism}) and we get
a morphism in $\mathcal{P}_{m}^{\alpha }$
\begin{equation*}
(\varphi _{\#}^{\ast },\varphi _{\#}):(R(\frak{G})^{\ast },R(\frak{G}%
))\rightarrow (R(\frak{H})^{\ast },R(\frak{H})).
\end{equation*}
By Theorem \ref{Gmod=} $_{\frak{H}}\mathcal{M}\approx \sigma \lbrack _{R(%
\frak{H})^{\ast }}R(\frak{H})],$ $_{\frak{G}}\mathcal{M}\approx \sigma
\lbrack _{R(\frak{G})^{\ast }}R(\frak{G})]$ and so we have the \emph{%
coinduction functor}
\begin{equation*}
\mathrm{Coind}_{\frak{H}}^{\frak{G}}(-):=\mathrm{Rat}^{R(\frak{G})}(\mathrm{%
Hom}_{R(\frak{H})^{\ast }-}(R(\frak{G})^{\ast },-):\text{ }{\normalsize _{%
\frak{H}}\mathcal{M}\rightarrow }\text{ }{\normalsize _{\frak{G}}}\mathcal{M}%
.
\end{equation*}
\end{punto}

\begin{lemma}
\label{jm}\emph{(\cite[Lemma 6.1.1, Corollary 6.1.2]{Swe69})} Let $%
I\vartriangleleft A$ be an ideal. If $_{A}I$ \emph{(}respectively $I_{A}$%
\emph{)} is finitely generated, then $_{A}I^{n}$ \emph{(}respectively $%
I_{A}^{n}$\emph{)} is finitely generated for every $n\geq 1.$ If moreover $%
I\subset A$ is $R$-cofinite, then $I^{n}\subset A$ is $R$-cofinite.
\end{lemma}

\begin{corollary}
\label{aff-hy}Let $\frak{G}$ be an affine monoid scheme \emph{(}respectively
an affine group scheme\emph{)} with coordinate ring $R(\frak{G}).$

\begin{enumerate}
\item  If $R$ is Noetherian, $_{R(\frak{G})}\omega $ is finitely generated
and $hy(\frak{G})\subset R^{R(\frak{G})}$ is pure, then $hy(\frak{G})$ is an
$R$-bialgebra \emph{(}respectively a Hopf $R$-algebra\emph{)} and $(R(\frak{G%
}),hy(\frak{G}))\in \mathcal{P}_{Big}^{\alpha }$ \emph{(}respectively $(R(%
\frak{G}),hy(\frak{G}))\in \mathcal{P}_{Hopf}^{\alpha }$\emph{)}.

\item  If $\frak{G}$ is \emph{infinitesimal flat}, then $hy(\frak{G})$ is an
\emph{infinitesimal flat} $R$-bialgebra \emph{(}respectively Hopf $R$-algebra%
\emph{)} and $(R(\frak{G}),hy(\frak{G}))\in \mathcal{P}_{Big}^{\alpha }$
\emph{(}respectively $(R(\frak{G}),hy(\frak{G}))\in \mathcal{P}%
_{Hopf}^{\alpha }$\emph{)}.
\end{enumerate}
\end{corollary}

\begin{Beweis}
\begin{enumerate}
\item  If $_{R(\frak{G})}\omega $ is finitely generated, then $\frak{F}%
_{\omega }\subset \mathcal{K}_{R(\frak{G})}$ by Lemma \ref{jm} and so $hy(%
\frak{G})\subset R(\frak{G})^{\circ }$ is an $R(\frak{G})$-subbimodule. The
result follows then from Proposition \ref{K(H)-zul} (1).

\item  The result follows from \cite[Lemma 9.2.1]{Mon93} and Proposition \ref
{K(H)-zul} (2).$\blacksquare $
\end{enumerate}
\end{Beweis}

\textbf{Acknowledgment.} This paper (up to a few changes) includes parts of
my doctoral thesis at the Heinrich-Heine Universit\"{a}t, D\"{u}sseldorf
(Germany). I am so grateful to my advisor Prof. Robert Wisbauer for his
wonderful supervision and the continuous encouragement and support. The
author thanks also the referee for his or her useful suggestions.

\end{document}